\documentclass[preprint,11pt]{elsarticle}

\newtheorem{theorem}{Theorem}[section]
\newtheorem{lemma}{Lemma}[section]
\newtheorem{corollary}{Corollary}[section]

\newtheorem{remark}{Remark}[section]
\newcommand{\ignore}[1]{}{}

\usepackage{amssymb}
\usepackage{amsmath}

\usepackage{color}
\usepackage{soul}
\usepackage[dvipsnames]{xcolor}

\usepackage[colorlinks=true, linkcolor=blue, citecolor=blue]{hyperref}

\def\1{{{\mbox{${\rm{1\negthinspace\negthinspace I}}$}}}}

\newcommand\beq{\begin{equation}}
\newcommand\eeq{\end{equation}}

\usepackage{ifthen}
\usepackage{xkeyval}
\usepackage{todonotes}
\begin{document}

\begin{frontmatter}

\title{Rates of convergence in the distances of Kolmogorov and Wasserstein for standardized martingales}
\author[cor1]{Xiequan Fan}
\author[cor2]{Zhonggen Su}
\address[cor1]{School of Mathematics and Statistics, Northeastern University at Qinhuangdao, Qinhuangdao, 066004, Hebei, China.}
\address[cor2]{School of Mathematical Sciences, Zhejiang University, Hangzhou 310058, China.}

\begin{abstract}
We give some rates of convergence in the distances of Kolmogorov and Wasserstein for standardized martingales
 with differences having finite variances. For the Kolmogorov distances, we present
 some exact Berry-Esseen bounds for  martingales,
 which generalizes some Berry-Esseen bounds due to Bolthausen.
For the Wasserstein distance, with Stein's method and Lindeberg's telescoping sum argument,  the rates of convergence
in martingale central limit theorems recover the classical rates for sums of i.i.d.\ random variables,
 and therefore they are believed to be optimal.

\end{abstract}

\begin{keyword} Martingales; Central limit theorem;  Rates of convergence;  Berry-Esseen bounds; Stein's method; the Wasserstein distance
\vspace{0.3cm}
\MSC Primary 60G42; 60F05;  Secondary 60E15
\end{keyword}

\end{frontmatter}




\section{Introduction}
\setcounter{equation}{0}

Let $\mathbf{X}=\{X_k, \mathcal{F}_k\}_{  k \geq 1}$ be a sequence of  square
integrable martingale differences defined on a probability space $(\Omega, \mathcal{F}, \mathbf{P})$, where $X_0=0$ and $\{\emptyset, \Omega\}=\mathcal{F}_0\subseteq ...\subseteq
\mathcal{F}_k\subseteq ...\subseteq \mathcal{F}$.
By the definition of martingale differences, we have
\begin{equation}
\mathbf{E}[X_k|\mathcal{F}_{k-1}]=0, \quad k \geq 1,  \ \ \textrm{a.s}.
\end{equation}
Define
\begin{equation}
S_0=0,\ \ \	S_k = X_1+X_2+...+X_k, \qquad k \geq 1.
\end{equation}
Then $\{S_k, \mathcal{F}_k\}_{ k \geq 1  }$ is a martingale.  Let $p \in [0,+\infty]$, and define
	\begin{equation}
		||\mathbf{X}||_p=\max_{ k \geq 1 } ||X_k||_p,
	\end{equation}
where $||X_k||_p=(\mathbf{E}|X_k|^p)^{1/p},    p \in [0,+\infty),$ and $||X_k||_{\infty}=\inf\{\alpha:\mathbf{P}(|X_k| \leq \alpha)=1\}$.
Let $n\geq 2.$ For briefness,  for  all $1\leq  k  \leq n,$ denote
	\begin{eqnarray*}
	&&	\sigma_k^2:=\mathbf{E}[X_k^2|\mathcal{F}_{k-1}], \ \ \ \  \ V_0^2:=0, \ \ \ \ \ \ \ V_k^2:=\sum_{i=1}^k \sigma_i^2, \ \ \ \ \ \ \ \rho_{n,k } :=\sqrt{V_n^2- V_{k-1}^2}\, ,\\
	 && \bar{\sigma}_k^2:=\mathbf{E} X_k^2 , \  \ \ \ \ \ \ \ \ \ \ \ \ \,   \ s_0^2:=0, \ \ \ \ \, \ \ \ \ s_k^2:=\sum_{i=1}^k \bar{\sigma}_i^2, \   \ \ \ \ \, \ \ \     \overline{\rho}_{n,k} :=\sqrt{s_n^2-s_{k-1}^2}\, ,  \\
	&& \upsilon_{n,k}:=\sqrt{\rho_{n,k}^2+a^2}  ,\ \ \ \ \ \ \ \   \,  \ \ \ \ \ \ \   \tau_{n,k} :=\sqrt{\overline{\rho}_{n,k}^2+a^2} .
	\end{eqnarray*}
Over all the paper, we assume that $s_n^2 \rightarrow \infty$ as  $n \rightarrow \infty$   and that, without loss of generality, $s_n^2 \geq 2$ for all $n$. The standardized martingale  is defined as $S_n/s_n .$
Denote by $\stackrel{\   \mathbf{P} \ }{\longrightarrow}$ convergence in probability.
The central limit theorem  (CLT) for martingales (see  the monograph Hall and Heyde \cite{HH80}) states that  the conditional ``Lindeberg condition",
that is for each $\varepsilon>0,$
\begin{eqnarray*}
  \frac{1}{s_n^2} \sum_{i=1 }^n \mathbf{E}[ X_i^2\mathbf{1}_{\{| X_i| \geq \varepsilon s_n  \}} | \mathcal{F}_{i-1}] \stackrel{\   \mathbf{P} \ }{\longrightarrow} 0, \ \ \ \ \   n \rightarrow \infty,
\end{eqnarray*}
and the ``conditional normalizing condition" $V_n^2/s_n^2   \stackrel{\   \mathbf{P}\   }{\longrightarrow} 1$,
together implies that the standardized martingales $ S_n/s_n$ converges in distribution to  the standard normal random variable.
In this paper, we are interested in the convergence rates in   the martingale CLT.

The distances of Kolmogorov and Wasserstein  are frequently used to describe  the convergence rates in   the martingale CLT.
Denote the Kologmorov  distance between $S_n/s_n$ and the standard normal random variable by
$$\textbf{K}(S_n/s_n)= \sup_{ x \in \mathbf{R}}\Big|\mathbf{P}(S_n/s_n \leq x)-\Phi \left( x\right) \Big| ,$$
where $\Phi \left( x\right)$ is the distribution function of the standard normal random variable.
Recall the definition of the $1$-Wasserstein distance.  The$1$-Wasserstein distance between two distributions $\mu$ and  $\nu$  is defined as follows:
\begin{eqnarray*}
	W_1(\mu, \nu)= \sup  \bigg \{ \Big|\mathbf{E}[f(X)] -\mathbf{E}[f(Y)]    \Big|:   \ (X, Y) \in \mathcal{L}(\mu, \nu),\ f\  \textrm{is  $1$-Lipschitz}     \bigg\} ,
\end{eqnarray*}
where  $\mathcal{L}(\mu, \nu)$ is the collection of all pairs of random variables whose marginal distributions are $\mu$ and $\nu$ respectively.
In particular, if $\mu_{S_n/s_n}$ is  the distribution of  $S_n/s_n$  and $\nu$ is the standard normal distribution, then
we write  $$ \textbf{W}\left(S_n/s_n\right) =W_1(\mu_{S_n/s_n} , \nu ). $$
The relationship between the distances of Kolmogorov and Wasserstein-1 can be found in Ross \cite{R11} (see p.\,218 therein).
 It is stated that
$$\textbf{K}(S_n/s_n) \leq \Big(  \frac{2}{  \pi } \Big)^{1/4}  \sqrt{\textbf{W}\left(S_n/s_n\right)} .$$
See also inequality (3.1) in Dedecker, \mbox{ Merlev\`{e}de}  and Rio \cite{DFR22}.

The exact convergence rates for $\textbf{K}(S_n/s_n)$, usually termed as Berry-Esseen's bounds, have been established
with various conditions.
For instance, when $p\in (1, 2]$, Heyde and Brown~\cite{HB70}~(see also Theorem~3.10 in Hall and Heyde~\cite{HH80}) have established the following Berry-Esseen bound: there exists a positive constant $C_p$ such that
\begin{equation}\label{sfdf}
 \mathbf{K}(S_n/s_n)  \leq  C_p \, N_n^{1/(2p+1) },
\end{equation}
where
\begin{equation}\nonumber
N_n = \sum_{i=1}^n \mathbf{E} |X_i/s_n|^{2p}  + \mathbf{E} \big| V^2_n/s_n^2 -1\big|^p.
\end{equation}
Later, Haeusler~\cite{H88} presented  an extension of (\ref{sfdf}), and proved that (\ref{sfdf}) holds  for all $p\in (1, \infty)$.
Moreover,   Haeusler \cite{H88} also gave an example  to show that his Berry-Esseen's bound is the best possible for martingales with finite moments.
When $ \mathbf{E} | V_n^2/s_n^2 -1|^{p} \rightarrow 0$ for some $p \in [1, \infty],$
the convergence rate  for $\mathbf{K}(S_n/s_n)$ has been established by Bolthausen \cite{2} (for $p=1, +\infty$),   Mourrat \cite{3},  El Machkouri and Ouchti \cite{EO07} and Fan~\cite{4}.
When $\sigma_k^2=\bar{\sigma}_k^2$ a.s.,   the  Berry-Esseen bound of Haeusler~\cite{H88} may not be the best possible.
At an earlier time, Bolthausen \cite{2} obtained the following Berry-Esseen bound. If there exist  three constants $\alpha,\ \beta$ and $\gamma$,
with $0 < \alpha \leq \beta < +\infty$ and $0 < \gamma < +\infty$,  such that
\begin{eqnarray}\label{sdfd}
\sigma_k^2=\bar{\sigma}_k^2\ \textrm{a.s.}, \ \ \ \  \alpha \leq \bar{\sigma}_k^2 \leq \beta \ \ \ \textrm{and} \ \ \   ||\mathbf{X}||_3 \leq \gamma,
\end{eqnarray}
 then  there exists a constant $C(\alpha,\beta,\gamma)$ depending only on $\alpha,\ \beta$ and  $\gamma$ such that
     \begin{eqnarray}\label{b}
		 \mathbf{K}(S_n/s_n) \leq C(\alpha,\beta,\gamma) \ n^{-1/4}.
	\end{eqnarray}
Bolthausen~\cite{2}~also proved the following result:
Assume that $V_n^2=s_n^2$ a.s.\ and that there exists a positive constant  $\gamma$ such that $||\mathbf{X}||_{\infty} \leq \gamma$.
Then there exists a constant $C(\gamma)$ depending on $\gamma$ such that
    \begin{eqnarray}\label{a}
		 \mathbf{K}(S_n/s_n) \leq C(\gamma) \frac{n \ln n}{s_n^3}.
	\end{eqnarray}
Moreover,  Bolthausen \cite{2} gave some examples to illustrate that the Berry-Esseen bounds (\ref{b}) and (\ref{a}) are sharp when  $s_n^2 $
is in order of $n$. The Berry-Esseen bound (\ref{a}) can be extended to a more general case.
Recently, Fan~\cite{4}~proved that if~$V_n^2=s_n^2$~a.s.\  and there exist two positive constants  $\delta$ and $\gamma$ such that
\begin{eqnarray}\label{BDer}
		\mathbf{E}[|X_k|^{2+\delta}|\mathcal{F}_{k-1}] \leq \gamma^{\delta} \mathbf{E}[ X_k ^2|\mathcal{F}_{k-1}]\ \ \ \textrm{a.s.},
\end{eqnarray}
then for any  $p \geq 1$, there exits a constant $C(p,\delta,\gamma)$ depending only on  $\delta,\ \gamma$ and $p$ such that
    \begin{eqnarray}\label{gfdgs1}
		 \mathbf{K}(S_n/s_n)\leq C(p,\delta,\gamma)\  \alpha_n  ,
    \end{eqnarray}
where
\begin{displaymath}
 \alpha_n = \left\{ \begin{array}{ll}
\displaystyle  s_n^{-\delta} , & \textrm{\ \ \ if $\delta \in (0,1)$,}\\
  & \textrm{ }\\
\displaystyle  s_n^{-1} \ln s_n , & \textrm{\ \ \ if $\delta \geq 1$.}
\end{array} \right.
\end{displaymath}
Fan \cite{4} also showed  that the last    Berry-Esseen bound is optimal under the stated condition.
When  $\delta=1$, an earlier result
  can be found in El Machkouri and Ouchti \cite{EO07}, where they obtained a Berry-Esseen bound of order $  s_n^{-1}  \ln n  $.
 When  the randomness for $\sigma_k^2$ tends to small as $k$ increasing, Bolthausen \cite{2} also gave the following
Berry-Esseen bound for standardized martingales.
Assume that   $ \sup_{k\geq 1}|| \mathbf{E}[ |X_k|^{3} | \mathcal{F}_{k-1} ] ||_\infty < \infty,$
and that there exists a positive constant $\sigma$ such that
 $\lim_{k \rightarrow \infty}\bar{\sigma}_k^2 =\sigma^2 $, and that there exists a positive  constant $\alpha\in (0, 1)  $ such that
 $$\mathbf{E} \big|\sigma_k^2-\overline{\sigma}_k^2 \big| = O( k^{-\alpha}).$$ Then
\begin{eqnarray}\label{ghd}
 \mathbf{K}(S_n/s_n) =O \bigg( \displaystyle   \frac{\ln n }{ n^{\alpha } } \bigg) .
\end{eqnarray}

Despite the fact that the convergence rates for $ \mathbf{K}(S_n/s_n)$ have been intensely studied,
there are only a few results for the convergence rates for $ \mathbf{W}(S_n/s_n)$.
To the best of our knowledge, we are only aware the work of Van Dung et al.\,\cite{5}, R\"ollin \cite{1}  and Dedecker,    Merlev\`{e}de  and Rio \cite{8,DFR22}. Van Dung et al.\,\cite{5} gave a generalization of  (\ref{a}) and proved that
(\ref{a}) holds also when $\mathbf{K}(S_n/s_n)$ is replaced by $\textbf{W}\big(S_{n} / s_{n}\big)$.
Using Stein's methods, R\"ollin \cite{1} proved the following result: If $V_{n}^{2}=s_n^2$ a.s.\ and $ \mathbf{E}\left|X_{i}\right|^{3} < \infty$ for all $1\leq i \leq n$, then for any $a\geq0$,
 \begin{eqnarray}\label{r}
 \textbf{W}\big(S_{n} / s_{n}\big) \leq \frac{3}{s_{n}}\sum_{i=1}^{n} \mathbf{E}\frac{\left|X_{i}\right|^{3}}{ \rho_{n,i }^{2}+a^2}+\frac{2a}{s_{n}}.
 \end{eqnarray}
 When $\mathbf{E} \big|\sigma_k^2-\overline{\sigma}_k^2 \big|\rightarrow 0,$ Dedecker,    Merlev\`{e}de  and Rio \cite{DFR22} recently also  have established a convergence rate for $ \textbf{W}(S_n/s_n),$
 with finite moment of order $2.$

In this paper, we are interested in giving  some generalizations of (\ref{b}),
(\ref{ghd})  and (\ref{r}) for martingales with finite moments of order  $2.$
From our results, we can recover many  optimal convergence rates for standardized martingales
with respect to the distances of Kolmogorov and Wasserstein.
Moreover, applications of our results are also discussed.
The paper is organized as follows. Our main results are stated and discussed in Section \ref{sec2}.
An application of our theorems to branching processes in a random environment is given in Section \ref{seca}. The
proofs of theorems and their corollaries are deferred to Section \ref{sec4}.

Throughout the paper,  $C$  will denote a finite and positive constant. This constant may vary from place to place.    For  two sequences of positive numbers $\{a_n\}_{n\geq 1}$ and $\{b_n\}_{n\geq 1}$,   write $a_n \asymp b_n$ if there exists a constant $C>0$ such that ${a_n}/{C}\leq b_n\leq C a_n$ for all sufficiently large $n$. We shall also use the notation $a_n \ll b_n$ to mean that there exists a positive constant $C$ not depending on $n$ such that $a_n \leq C b_n$.
We also  write $a_n \sim b_n$ if  $\lim_{n\rightarrow \infty} a_n/b_n =1$.
We denote by $N(0, \sigma^2)$   the normal  distribution with mean $0$ and variance $\sigma^2$.
$\mathcal{N}$  will designate a standard normal random variable, and we will denote by $\Phi(\cdot)$  the cumulative distribution function of a standard normal random variable.

\section{Main results}\label{sec2}
\setcounter{equation}{0}
We have the following convergence rates for standardized martingales $S_n/s_n$ in the distances of Kolmogorov  and Wasserstein.

\subsection{Convergence rates in the Kolmogorov   distance}
We have the following Berry-Esseen bound for the standardized martingale $S_n/s_n$.
\begin{theorem}\label{th4.s1}
Assume that $|| \textit{\textbf{X}} ||_2 < \infty$.
 Let $a\geq 1$.   Then it holds
\begin{eqnarray*}
 \mathbf{K}(S_n/s_n)
   \  \ll \  \sum_{k=1}^n \Bigg(  \frac{ \mathbf{E}\big[  |X_k|^3  \mathbf{1}_{\{|X_k|\leq \tau_{n,k+1}\}} \big]   }{\tau_{n,k+1}^3}+  \frac{ \mathbf{E}\big[  X_k ^2  \mathbf{1}_{\{|X_k|> \tau_{n,k+1}\}} \big] + \mathbf{E}   \big | \sigma_k^2   - \overline{\sigma}_k^2       \big| }{\tau_{n,k+1}^2}   \Bigg) +    \frac{a}{s_n }.
\end{eqnarray*}
\end{theorem}

By Theorem \ref{th4.s1},   the following CLT for the standardized martingales $ S_n/s_n$ holds. If
there exists a constant $a\geq 1$ such that
\begin{eqnarray}\label{gsdxvn}
 \sum_{k=1}^n \Bigg(  \frac{ \mathbf{E}\big[  |X_k|^3  \mathbf{1}_{\{|X_k|\leq \tau_{n,k+1}\}} \big]   }{\tau_{n,k+1}^3}+  \frac{ \mathbf{E}\big[  X_k ^2  \mathbf{1}_{\{|X_k|> \tau_{n,k+1}\}} \big] }{\tau_{n,k+1}^2}   \Bigg) \rightarrow 0 \ \ \ \textrm{and} \ \ \  \  \sum_{k=1}^n   \frac{ \mathbf{E}   \big | \sigma_k^2   - \overline{\sigma}_k^2       \big| }{\tau_{n,k+1}^2} \rightarrow 0,
\end{eqnarray}
then  $ S_n/s_n$ converges in distribution to the standard normal distribution as $n\rightarrow \infty$.

For $|| \textit{\textbf{X}} ||_p < \infty,  p\in (2, 3],$ a closely related result to Theorem \ref{th4.s1} can be found in Dedecker,  Merlev\`{e}de and Rio \cite{DFR22}, see Corollary 3.1 therein.
The next corollary  indicates that in certain cases,  the Berry-Esseen bound  in Theorem \ref{th4.s1} is sharper than the one in Corollary 3.1 of Dedecker,  Merlev\`{e}de and Rio \cite{DFR22}.

\begin{corollary}\label{co2.5}
Assume that $|| \textit{\textbf{X}} ||_{2 } < \infty$, and that there exists a constant $\delta\in (0, 1]$ such that
$$\mathbf{E}  |X_k|^{2+\delta}  \leq C \,\mathbf{E} X_k^2,\ \ \ 1\leq k \leq n.$$
 If there exists a positive constant $\alpha  $ such that
 $$\mathbf{E} \big|\sigma_k^2-\overline{\sigma}_k^2 \big| \leq C \, \frac{ \overline{\sigma}_k^2}{\, s_n^{\alpha} \, } ,\ \ \ 1\leq k \leq n,$$ then  it holds
\begin{displaymath}
  \mathbf{K}(S_n/s_n) \ll  \,  \displaystyle \frac{ 1}{ s_n^{\delta /(1+ \delta) }}+ \frac{\ln s_n  }{ s_n^{\alpha } }    .
\end{displaymath}
\end{corollary}

When $\sigma_k^2=\overline{\sigma}_k^2$ a.s., $\delta=1$ and $\mathbf{E}  |X_k|^{3}  \leq C \,\mathbf{E} X_k^2$ ($1\leq k \leq n$),
Dedecker,  Merlev\`{e}de and Rio \cite{DFR22}  obtained a Berry-Esseen bound of order $ s_n^{-1/2} \sqrt{\ln s_n}$, see the remark after Corollary 3.1 in \cite{DFR22}. Now, Corollary \ref{co2.5} gives  a Berry-Esseen bound of order $ s_n^{-1/2},$ by (\ref{b}), which is sharp when  $s_n^2\asymp n$.

By Theorem \ref{th4.s1}, we can obtain the following corollary, which generalizes inequality (\ref{b}).
\begin{corollary}  \label{th4.21}
Assume that there exists a constant $ \delta \in (0, 1]$ such that $||\mathbf{X}||_{2+\delta} < \infty$,  and  that
 $\inf_{k \geq 1}\bar{\sigma}_k^2  \geq\sigma^2$ a.s., with $\sigma$ a positive constant.
If there exists a positive constant $\alpha \in (0, 1) $ such that
$$\mathbf{E} \big|\sigma_k^2-\overline{\sigma}_k^2 \big| \ll k^{-\alpha }   ,\ $$ then  it holds for all $n\geq 2,$
\begin{displaymath}
  \mathbf{K}(S_n/s_n) \ \ll \ \displaystyle \frac{ 1}{ n^{\delta /(2+2\delta) }}+ \frac{\ln  n  }{ n^{\alpha } }  .
\end{displaymath}
\end{corollary}

Clearly, when $\delta=1$ and $\sigma_k^2=\overline{\sigma}_k^2$ a.s., the last corollary reduces to the Berry-Esseen bound
(\ref{b}).  Thus, Corollary \ref{th4.21} can be regarded as a generalization of
 (\ref{b}), that is Theorem 1 of Bolthausen~\cite{2}.

For $\delta>0$ and $n \geq 20$,   Dedecker,  Merlev\`{e}de and Rio \cite{DFR22} showed that there exists a  finite sequence of martingale differences $(\Delta M_i, \mathcal{F}_i)_{1 \leq i  \leq n}$ satisfying the following three conditions:
\begin{enumerate}
  \item  $\mathbf{E} [ \Delta M_i^2 | \mathcal{F}_{i-1} ] =1$ a.s.,
  \item  $\sup_{1 \leq i \leq n}\mathbf{E}   |\Delta M_i|^{2+\delta}     \leq  \mathbf{E} | \mathcal{N}|^{2+\delta}   + 5^{\delta}$,
  \item   $\sup_{u \in \mathbf{R}} \Big| \mathbf{P}\Big(    \sum_{i=1}^n \Delta M_i \leq u  \sqrt{n}\Big) -\Phi (x)\Big| \geq  0.06 \, n^{-\delta/(2+2\delta)}.  $
\end{enumerate}
Under the conditions  1 and 2, by Corollary \ref{th4.21} with $\alpha=1/2$,
we have for all $\delta\in (0, 1],$
\begin{eqnarray*}
 \sup_{u \in \mathbf{R} }  \bigg|  \mathbf{P}\Big(    \sum_{i=1}^n \Delta M_i \leq u  \sqrt{n}\Big)  -  \Phi (u)\bigg|
 \  \ll \  \frac{ 1  }{ n^{\delta /(2+2\delta) }} .
\end{eqnarray*}
By point 3 and the last inequality, we find that
the decaying rate $ n^{-\delta /(2+2\delta) } $ in the Berry-Esseen bound of Corollary \ref{th4.21}  is sharp.

\begin{theorem}\label{th4.2}
Assume  $  \sup_{k\geq 1}|| \mathbf{E}[  X_k ^2 | \mathcal{F}_{k-1} ] ||_\infty < \infty$,   $\rho_{n,k}^2 \leq  C\, \overline{\rho}^2_{n,k}$ for all $ k $ and
\begin{eqnarray} \label{fsdfzn}
 \lim_{a \rightarrow \infty} \sum_{k=1}^n\Bigg(  \frac{ \big\| \mathbf{E} [  |X_k|^3  \mathbf{1}_{\{|X_k|\leq \tau_{n,k+1}\}} |\mathcal{F}_{k-1} ]  \big\|_{\infty}}{\tau_{n,k+1}^3}  +  \frac{ \big\| \mathbf{E} [  X_k^2  \mathbf{1}_{\{|X_k|> \tau_{n,k+1}\}} |\mathcal{F}_{k-1} ]  \big\|_{\infty}}{\tau_{n,k+1}^2}  \Bigg)  =0
\end{eqnarray}
uniformly for $n \geq 2.$
Then it holds for all $a$ large enough,
\begin{eqnarray*}
 \mathbf{K}(S_n/s_n)
     \!&\ll&\!   \frac{  1 }{s_n} \sum_{k=1}^n \Bigg(  \frac{ \big\| \mathbf{E} [  |X_k|^3  \mathbf{1}_{\{|X_k|\leq \tau_{n,k+1}\}} |\mathcal{F}_{k-1} ]  \big\|_{\infty} }{\tau_{n,k+1}^2 } \, + \,  \frac{ \big\|\mathbf{E} [  X_k ^2  \mathbf{1}_{\{|X_k|> \tau_{n,k+1}\}} |\mathcal{F}_{k-1} ]\big\|_{\infty}   }{\tau_{n,k+1} } \Bigg)  \nonumber  \\
  && \   + \   \sum_{k=1}^n \frac{   \mathbf{E}   \big | \sigma_k^2   - \overline{\sigma}_k^2       \big| }{\tau_{n,k+1}^2 }    +    \frac{a}{s_n } .
\end{eqnarray*}
\end{theorem}

The condition (\ref{fsdfzn}) can be obtained if, for instance, $ || \mathbf{E}[ |X_k|^{2+\delta} | \mathcal{F}_{k-1} ] ||_\infty  \leq C \, \overline{\sigma}_k^2$, $ \delta>0$.
Thus, we have the following corollary.

\begin{corollary}\label{fssdvs}
Assume that  $||\mathbf{X}||_2 < \infty$, and that there exists a positive constant $ \delta  $ such that
  $$|| \mathbf{E}[ |X_k|^{2+\delta} | \mathcal{F}_{k-1} ] ||_\infty  \leq C \, \overline{\sigma}_k^2$$  and   $\rho_{n,k}^2 \leq  C\, \overline{\rho}^2_{n,k}$ for all $k$. If there exists a positive  constant $\alpha $ such that for all $1\leq k \leq n,$
 $$\mathbf{E} \big|\sigma_k^2-\overline{\sigma}_k^2 \big| \leq C \, \frac{ \overline{\sigma}_k^2}{\, s_n^{\alpha} \, } ,$$ then the
 following   inequalities  hold.
 \begin{itemize}
\item[\emph{(i)}] If~$\delta \in (0, 1)$, then
     \begin{eqnarray}
	  \mathbf{K}(S_n/s_n)  \ll   \displaystyle \frac{ 1}{ s_n^{\delta  }}+ \frac{\ln s_n  }{ s_n^{\alpha } }    .
	\end{eqnarray}
\item[\emph{(ii)}] If~$\delta\geq 1$, then
     \begin{eqnarray} \label{rdffs}
	  \mathbf{K}(S_n/s_n)   \ll  \displaystyle \frac{ \ln s_n  }{ s_n }+ \frac{\ln s_n  }{ s_n^{\alpha } }    .
	\end{eqnarray}
\end{itemize}
\end{corollary}

When $s_n^2$ is in order of $n,$ the last corollary implies the following result, which shows that our results are optimal.
\begin{corollary}  \label{fsdvs}
Assume that
 $\inf_{k \geq 1}\bar{\sigma}_k^2  >0$,  and that there exists a positive constant $ \delta  $ such that $$ \sup_{k\geq 1}|| \mathbf{E}[ |X_k|^{2+\delta} | \mathcal{F}_{k-1} ] ||_\infty < \infty.$$
 If there exists a positive  constant $\alpha\in (0, 1)  $ such that
 $$\mathbf{E} \big|\sigma_k^2-\overline{\sigma}_k^2 \big| = O( k^{-\alpha}),$$ then the
 following two inequalities  hold.
\begin{itemize}
\item[\emph{(i)}] If $\delta \in (0, 1)$, then it holds
\begin{eqnarray}\label{rdfdf}
\mathbf{K}(S_n/s_n) \ll \displaystyle \frac{ 1}{ n^{\delta/2  }}+ \frac{\ln n  }{  n^{\alpha } }   .
\end{eqnarray}

\item[\emph{(ii)}] If  $\delta \geq 1$, then it holds
\begin{eqnarray} \label{rdfdfs}
 \mathbf{K}(S_n/s_n)  \ll \displaystyle \frac{ \ln n  }{  \sqrt{n}\ }+ \frac{\ln n }{ n^{\alpha } }  .
\end{eqnarray}
\end{itemize}
\end{corollary}

When $\delta=1,$  Bolthausen \cite{2} obtained a result similar to (\ref{rdfdfs}), cf.\ Theorem 3~(a) therein.
As Corollary \ref{fsdvs} also includes the case $\delta \in (0, 1),$  Corollary \ref{fsdvs} can be
regarded as a generalization of Theorem 3~(a) of Bolthausen \cite{2}.

The convergence rates $ \frac{ 1}{ n^{\delta/2  }}$ and $  \frac{ \ln n  }{  \sqrt{n} }$ in the inequalities  (\ref{rdfdf}) and  (\ref{rdfdfs}), respectively, are the same as that for independent and identically distributed (i.i.d.) random variables, up to a factor $\ln n$. However, the factor $\ln n$ cannot be replaced by $1
$. Indeed, by Example 2 of Bolthausen \cite{2}, the convergence rate $\frac{\ln n  }{  \sqrt{n}}$ in (\ref{rdfdfs}) is the best possible even for
uniformly bounded martingale differences with $\sigma_k^2=\overline{\sigma}_k^2=1$ a.s.
Thus the terms $\frac{ 1}{ n^{\delta/2  }}$ and $\frac{ \ln n  }{  \sqrt{n} }$ in the inequalities  (\ref{rdfdf}) and  (\ref{rdfdfs})  are optimal.

It is known that there exists a martingale such that $\|\sigma_k^2-\overline{\sigma}_k^2 \|_{\infty} = O( k^{-\alpha})$ and
$\displaystyle \limsup _{n \rightarrow \infty} n^{\alpha } \mathbf{K}(S_n/s_n) > 0,$ see Example 4 of Bolthausen \cite{2}. Thus, the
the terms  $ \frac{\ln n }{ n^{\alpha } }$ in  Corollary  \ref{th4.21}, the inequalities  (\ref{rdfdf}) and  (\ref{rdfdfs})  are optimal, up to a factor $\ln n
$.

\begin{remark}
If the condition $\inf_{k \geq 1}\bar{\sigma}_k^2  >0$ is replaced by $  \liminf_{k\rightarrow \infty}\bar{\sigma}_k^2  =\sigma^2$, with $\sigma$ a positive constant,   then  Corollaries \ref{th4.21} and \ref{fsdvs} hold also. The proofs are similar.
\end{remark}

\subsection{Convergence rates in the Wasserstein  distance}
When $V_n^2=s_n^2$ for some $n$,
using Stein's method  (see Chen, Goldstein and Shao \cite{7}) and Lindeberg's telescoping sum argument (see Bolthausen \cite{2}), we obtain the following convergence rates in the Wasserstein  distance for standardized martingales  $S_n/s_n$.
\begin{theorem}\label{theorem2.1}
Assume that $V_n^2=s_n^2$ for some $n$, and $||\mathbf{X}||_2 < \infty$. Then for any $a \geq 0$, it holds
\begin{eqnarray*}
 && \textbf{W}\big(  S_n/s_n  \big)  \leq\frac{1}{s_n} \sum_{k=1}^n  \Bigg(\mathbf{E}\bigg[\frac{(\sigma_k^2  + X_k^2)\mathbf{1}_{\{|X_k| \geq \upsilon_{n,k} \}} }{\upsilon_{n,k}} \bigg]   +2\, \mathbf{E}\bigg[ \frac{ (\sigma_k^2 |X_k|+|X_k|^3)\mathbf{1}_{\{|X_k| <\upsilon_{n,k}\}}}{\upsilon_{n,k}^2}\bigg]\Bigg)      +\frac{2a}{s_n}.\ \ \
\end{eqnarray*}
\end{theorem}

Notice that in the right hand side of the last inequality, the constants are explicit. This is an advantage of Stein's method.

When the martingale differences have finite moments of order $2+\delta, \delta \in (0, 1], $
Theorem \ref{theorem2.1} implies the following corollary.
\begin{corollary} \label{fdgs}
Assume that $V_n^2=s_n^2$ a.s.\ for some $n$, and that there exists a constant $\delta \in (0, 1]$ such that
$||\mathbf{X}||_{2+\delta} < \infty$. Then for any $a \geq 0$, it holds
\begin{eqnarray*}
 \textbf{W} \big( S_n/s_n  \big)   \leq \frac{6}{s_n} \sum_{k=1}^n    \mathbf{E} \frac{ |X_k|^{2+\delta} }{  (\rho_k^2+a^2 )^{(1+ \delta)/2} }        +\frac{ 2a}{s_n}.
\end{eqnarray*}
\end{corollary}

Clearly, when $\delta=1,$   Corollary  \ref{fdgs} reduces to the result of  R\"ollin \cite{1},
up to a constant.

When the conditional moments of order $2+\delta$ can be dominated by the conditional variances,
Theorem  \ref{theorem2.1} implies the following corollary, which extends the main result of Fan \cite{4}
 to the distance of Wasserstein.

\begin{corollary} \label{fdfff}
Assume that $V_n^2=s_n^2$ a.s.\ for some $n \geq 2$, and that there exists a positive constant $\gamma$
such that for all $1\leq k \leq n,$
$$\mathbf{E}[ |X_k|^{2+\delta} | \mathcal{F}_{k-1} ] \leq \gamma \, \mathbf{E}[ X_k^2 | \mathcal{F}_{k-1} ].$$
Then we have the following two inequalities.
\begin{itemize}
\item[\emph{(i)}] If $\delta \in (0, 1)$, then
     \begin{eqnarray}
	 \textbf{W} \big(  S_n/s_n \big)  \leq C \, s_n^{-\delta } .
	\end{eqnarray}

\item[\emph{(ii)}] If $\delta \geq 1$, then
     \begin{eqnarray}\label{rffds}
	 \textbf{W} \big( S_n/s_n  \big)   \leq C \,  \frac{\ln s_n  }{ s_n } .
	\end{eqnarray}
\end{itemize}
\end{corollary}

  R\"ollin \cite{1} (cf.\ Corollary 2.3 therein) considered the case   $\mathbf{E}[ |X_k|^3 | \mathcal{F}_{k-1} ] \leq \gamma \, \mathbf{E}[ X_k^2 | \mathcal{F}_{k-1} ]$  and obtained a convergence rate $\displaystyle C\,  \frac{  \ln n}{ s_n }$.
Because the condition $\mathbf{E}[ |X_k|^3 | \mathcal{F}_{k-1} ] \leq \gamma \, \mathbf{E}[ X_k^2 | \mathcal{F}_{k-1} ]$ implies  $s_n^2 \leq n \gamma^2,$ we get $\displaystyle  \frac{  \ln s_n}{ s_n } =O \Big(\frac{ \ln n}{ s_n } \Big )$. Thus  (\ref{rffds}) slightly
improved the result in Corollary 2.3 of \cite{1}.
Moreover, Corollary \ref{fdfff} can be regarded as a generalization of Corollary 2.3 of R\"ollin \cite{1} to the case $\delta \in (0, 1)$.

When the martingale differences satisfy the conditions for the Berry-Esseen bound (\ref{b}) in Bolthausen \cite{2},  Corollary \ref{fdgs} implies the following result.
\begin{corollary}\label{fdffs}
Assume that $V_n^2=s_n^2$ a.s.\ for some $n \geq 2$ and that there exists two positive constants
$\alpha$  and $\delta$   such that  $\sigma_k^2 \geq \alpha$ a.s.\ for all $1\leq k \leq n$ and  $||\mathbf{X}||_{2+\delta} < \infty.$ Then the following two inequalities hold.
\begin{itemize}
  \item[\emph{(i)}] If $\delta \in (0, 1)$, then
     \begin{eqnarray}
	 \textbf{W} \big(  S_n/s_n  \big) \ll   \frac{1}{ n^{\delta/2}}.
	\end{eqnarray}

  \item[\emph{(ii)}] If $\delta \geq 1$, then
     \begin{eqnarray}\label{rffs}
	 \textbf{W} \big(  S_n/s_n  \big)   \ll    \frac{\ln n  }{ \sqrt{n}\ } .
	\end{eqnarray}
\end{itemize}
\end{corollary}

The convergence rates in Corollary \ref{fdffs} are the same as that for i.i.d.\ random variables, up to the factor $\ln n$ in the bound (\ref{rffs}). Thus, up to the factor $\ln n$ for the case $\delta =1$, the convergence rates in Corollary \ref{fdffs} are optimal.

Comparing the results of  (\ref{b}) and Corollary \ref{fdffs}, we find an interesting conclusion: when $\delta=1$,
the convergence rate for standardized martingales in the Wasserstein distance is much better than the one in the Kolmogorov  distance.

When $\mathbf{E} \big|\sigma_k^2-\overline{\sigma}_k^2 \big|\rightarrow 0\, (k\rightarrow \infty)$, that is the randomness for $\sigma_k^2$ tends to small as $k$ increasing,
Dedecker, Merlev\`{e}de and Rio \cite{DFR22}  have  established a convergence rate  in the Wasserstein distance for standardized martingales  $S_n/s_n$. With Stein's method, we can establish a result similar to that of Dedecker, Merlev\`{e}de and Rio \cite{DFR22}, but with   explicit constants.
We have the following result.
\begin{theorem}\label{theorem3.2}
Assume that $||\mathbf{X}||_2 < \infty$. For any $a \geq 0$, it holds
\begin{eqnarray*}
  \textbf{W} \big(  S_n/s_n  \big)
	 \!\! &\leq&\! \!  \frac{1}{s_n} \sum_{k=1}^n \Bigg(\frac{ \mathbf{E} \big|\sigma_k^2-\overline{\sigma}_k^2 \big| +  \mathbf{E}[ X_k^2\mathbf{1}_{\{|X_k| \geq \tau_{n,k}\}}] }{\tau_{n,k}}\,   \\
 &&\ \ \ \ \ \ \ \ \ \ \ \  + \, \frac{   3\,\overline{\sigma}_k^2  \mathbf{E}  |X_k| + 2\,\mathbf{E}[|X_k|^3\mathbf{1}_{\{|X_k| <\tau_{n,k}\}}  ]}{\tau_{n,k}^2 } \Bigg) +\frac{2a}{s_n}.
\end{eqnarray*}
\end{theorem}

Notice that $\lim_{n\rightarrow \infty} \sup_{k}\overline{\sigma}_k^2/s_n^2 =0$. As $s_n\geq2,$  we have for all  $ a \geq 1,$ $$\displaystyle \sum_{k=1}^n\frac{  \overline{\sigma}_k^2}{ \tau_{n,k}^2  } \mathbf{E}  |X_k| \leq \sum_{k=1}^n\frac{  \overline{\sigma}_k^2/s_n^2}{ \overline{\rho}_{n,k}^2/s_n^2 +a^2/s_n^2 }  \sqrt{ ||\mathbf{X}||_2}  \leq C   \ln   s_n    .$$
Applying the last line to Theorem \ref{theorem3.2} with $a \geq 1,$ we can deduce that as $s_n\rightarrow \infty,$
\begin{eqnarray}
  \textbf{W} \big(  S_n/s_n  \big)
	 \ll   \frac{1}{s_n} \sum_{k=1}^n \Bigg(\frac{ \mathbf{E} \big|\sigma_k^2-\overline{\sigma}_k^2 \big|+  \mathbf{E}\big[ X_k^2\mathbf{1}_{\{|X_k| \geq \tau_{n,k}   \}} \big] }{\tau_{n,k}    }  \,  + \, \frac{ \mathbf{E}\big[|X_k|^3\mathbf{1}_{\{|X_k| <\tau_{n,k}   \}} \big ]  }{ \tau_{n,k}^2  }\Bigg)   +   \, \frac{\ln  s_n}{s_n}  . \label{gsds}
\end{eqnarray}
From the last inequality, we get the following convergence rate.
\begin{corollary} \label{fsvg}
Assume that $||\mathbf{X}||_2 < \infty$ and that  $\inf_{k \geq 1}\bar{\sigma}_k^2   \geq\sigma^2$, with $\sigma$ a positive constant.
If there exists a constant $\alpha \in (0, 1]$ such that $$\mathbf{E} \big|\sigma_k^2-\overline{\sigma}_k^2 \big| \ll k^{-\alpha},$$ then for any $a \geq 1$, it holds
\begin{eqnarray*}
 \textbf{W} \big(  S_n/s_n  \big)   \ll    \frac{1}{ \sqrt{n}} \sum_{k=1}^n \Bigg(\frac{  \mathbf{E}[ X_k^2\mathbf{1}_{\{|X_k| \geq\tau_{n,k}\}}] }{\tau_{n,k} }  \,    + \, \frac{ \mathbf{E}[|X_k|^3\mathbf{1}_{\{|X_k| <\tau_{n,k}\}}  ]  }{ \tau_{n,k}^2 }\Bigg) +\frac{\ln n}{ n^{\alpha} } +  \frac{\ln  n}{ \sqrt{n}\, }  .
\end{eqnarray*}
Moreover, if $||\mathbf{X}||_{2+\delta} < \infty$,  then the following two inequalities hold.
\begin{itemize}
\item[\emph{(i)}] If $\delta \in (0, 1)$, then
\begin{eqnarray}
\textbf{W} \big(  S_n/s_n  \big) \ll    \frac{1}{ n^{\delta/2 }} + \frac{\ln n}{ n^{\alpha} }   .
\end{eqnarray}
\item[\emph{(ii)}] If  $\delta \geq 1$, then
\begin{eqnarray} \label{rdffs}
\textbf{W} \big( S_n/s_n  \big)   \ll  \frac{\ln n  }{ \sqrt{n}\,  } + \frac{\ln n}{ n^{\alpha} }   .
\end{eqnarray}
\end{itemize}
\end{corollary}

The convergence rates $\displaystyle  \frac{1}{ n^{\delta/2 }}$  and $\displaystyle  \frac{\ln n  }{ \sqrt{n} \, }$ in Corollary \ref{fsvg} are the same as the classical rates for i.i.d.\ random variables, up to the factor  $\ln n$. Thus, the convergence rates in Corollary \ref{fsvg} are
optimal or almost optimal.

\section{Application to branching processes in a random environment} \label{seca}
 \setcounter{equation}{0}
The branching processes  in a random environment can be described as follows.
Let $\xi=\{\xi_n\}_{n\geq0} $ be a sequence of i.i.d.\ random variables, where $\xi_n$ stands for the random environment at generation $n.$  For each $n \in \mathbf{N}=\{0, 1,...\}$,  the realization of $\xi_n$ corresponds a probability law $\{ p_i(\xi_n): i \in  \mathbf{N}\}.$
 A branching process $\{Z_n\}_{n\geq 0}$ in the random environment $\xi$ can be defined as follows:
 \begin{equation}
 Z_0=1,\ \ \ \ Z_{n+1}= \sum_{i=1}^{Z_n} X_{n,i}, \   \ \ \ \  n \geq 0,
 \end{equation}
where $X_{n,i}$ is the offspring number of   the $i$-th individual in the generation $n.$
Moreover, given the environment $\xi,$ the random variables $  \{X_{n,i}\}_{i\geq 1} $ are independent of each other with a common  distribution law
as follows:
  \begin{equation} \label{ssfs}
 \mathbf{P}(X_{n,i} =k | \xi   )     = p_k(\xi_n),\ \ \ \ \ \ \ k\in \mathbf{N},
\end{equation}
and they are also independent of $\{Z_i\}_{1\leq i \leq n}.$
Denote by $\mathbf{P}_\xi$ the quenched law, i.e.\ the conditional probability when the environment
$\xi$ is given, and by $\tau$  the law of the environment $\xi.$ Then $\mathbf{P}(dx, d\xi )=\mathbf{P}_\xi(dx)\tau(d\xi)$
is the total law of the process $\{Z_n\}_{n\geq0}$, called annealed law. In the sequel, the expectations with respect to the quenched law and annealed  law are denoted
by $\mathbf{E}_\xi$ and $\mathbf{E}$, respectively.

  Denote by $m$  the average offspring number  of an individual, usually termed the offspring mean.
 Clearly, it holds  $$m=\mathbf{E}Z_1=\mathbf{E} X_{n,i} =\sum_{k=1}^\infty  k \mathbf{E} p_k(\xi_{n,i}).$$
Thus, $m$ is not only the mean of $Z_1$, but also the mean of each individual $X_{n,i}$.
Denote
\begin{eqnarray*}
  m_n:=\mathbf{E}_\xi X_{n,i}=\sum_{i=1}^{\infty} i \, p_i(\xi_n) ,\ \   \ n\geq 0,   \ \ \ \ \ \mu= \mathbf{E} \ln m_0,  \\
\ \ \ \ \nu^2= \mathbf{E}(  \ln m_0 -\mu)^2 \ \ \ \ \ \ \ \ \ \ \  \   \textrm{and} \ \ \ \ \ \ \ \ \ \   \ \ \ \  \tau^2=\mathbf{E} ( Z_1-m_0  )^2.
\end{eqnarray*}
Notice that  $\{m_n\}_{n\geq 1}$ is  a sequence of nonnegative i.i.d.\ random variables. The distribution of $m_n$ depends only on $\xi_n$, and
$m_n$ is independent of $\{Z_i\}_{1\leq i \leq n}$. Clearly, it holds $m= \mathbf{E}m_n.$ The process $\{Z_n \}_{ n\geq 0}$  is respectively  called supercritical, critical or subcritical according to $\mu >0$, $\mu =0$ or $\mu <0$. 

Assume
\begin{eqnarray}\label{defv1}
 p_0(\xi_0)=0  \ \ \ \ a.s.,
\end{eqnarray}
which means each individual has at least one offspring.
Thus the process $\{Z_n\}_{n\geq 0 }$ is   supercritical.
 Denote by $v$  and $\sigma$ the standard variances of $Z_1$ and $m_0$ respectively,  that is
\begin{eqnarray*}
  \upsilon^2=\mathbf{E} (Z_1-m)^2,\ \ \ \ \ \ \sigma^2= \mathbf{E} (m_0-m)^2.
\end{eqnarray*}
To avoid triviality,   assume that  $0 < v, \sigma  <\infty,$ which implies that
  $\mu, \nu^2$ and $\tau$ all are finite.
The assumption $\sigma  >0$  means that the random environment is not degenerate.
For the degenerate case $\sigma =0$ (i.e.\ Bienaym\'{e}-Galton-Watson process),
Maaouia and Touati \cite{MT05} have established  the  CLT for the maximum likelihood estimator of $m$.
Moreover, we also assume that
\begin{eqnarray}\label{defv3}
\mathbf{E} \frac{Z_1}{m_0} \ln^+ Z_1  < \infty.
\end{eqnarray}
Denote $V_n= \frac{Z_n }{\Pi_{i=0}^{n-1}m_i},\   n \geq 1. $
The  conditions (\ref{defv1}) and  (\ref{defv3})  imply
that the limit $V=\lim_{n\rightarrow \infty} V_n $ exists a.s.,
$V_n \rightarrow V$ in $\mathbf{L}^{1} $ and   $\mathbf{P}(V >0)=\mathbf{P}(\lim_{n \rightarrow \infty} Z_n = \infty)=\lim_{n \rightarrow \infty} \mathbf{P}(Z_n >0) =1$
(see   Athreya and Karlin \cite{AK71b} and   Tanny \cite{T88}).
We also need the following assumption: There exist  two constants $p > 1$ and $\eta_0 \in (0, 1)$ such that
\begin{eqnarray}\label{defv4}
\mathbf{E} (\theta_0(p))^{\eta_0} < +\infty, \ \ \ \ \ \
\textrm{where} \ \ \  \theta_0(p)= \frac{Z_1^p}{m_0^p} .
\end{eqnarray}
With the last assumption, Grama, Liu  and Miqueu \cite{GLP20} proved that there exists a positive constant $\alpha$ such that
\begin{eqnarray}\label{defv55}
\mathbf{E} V^{-\alpha}   < \infty.
\end{eqnarray}

A  critical task in statistical inference of BPRE is to estimate  the offspring mean   $m.$
To this end,   the Lotka-Nagaev    \cite{L39,N67}   estimator $\frac{Z_{ k+1}}{Z_{ k}}$  plays an important role.
  Denote the average  Lotka-Nagaev  estimator   by $\hat{m}_{n_0,n}$, that is
  $$\hat{m}_{n_0,n}=\frac{1}{n}\sum_{k=n_0}^{n_0+n-1} \frac{Z_{ k+1}}{Z_{ k}}.$$
 For any $n\geq2$ and $n_0 \geq 0$,  denote
\begin{eqnarray}\label{fdsds}
S_{n_0,n}  = \frac{ \sqrt{n}}{ \sigma   \ }  (\hat{m}_{n_0,n} -m )
\end{eqnarray}
the normalized process for $\hat{m}_{n_0,n}$. Here, $n_0$ may depends on $n.$
We have the following  Berry-Esseen bound for $S_{n_0,n}$ and $-S_{n_0,n}$.
\begin{theorem}\label{th00s1}
Assume that there exists   a positive constant $\delta  $ such that
\begin{equation} \nonumber
\mathbf{E}| m_0 -m   |^{2+\delta } + \mathbf{E}| Z_1  -m_0   |^{2+\delta }   < \infty.
\end{equation}
\begin{itemize}
\item[\emph{(i)}] If~$\delta \in (0, 1)$, then
     \begin{eqnarray}\label{rdffs1}
	 \mathbf{K}(  S_{n_0,n}) \ll   \displaystyle \frac{ 1}{ n^{\delta/2}} \ \ \    \ \ \ \ \ and \ \ \    \ \ \ \ \textbf{W} (  S_{n_0,n})  \ll   \displaystyle \frac{ 1}{ n^{\delta/2}}.
	\end{eqnarray}
\item[\emph{(i)}] If~$\delta\geq 1$, then
     \begin{eqnarray} \label{rdffs2}
	 \mathbf{K}(  S_{n_0,n})   \ll  \displaystyle \frac{ \ln n  }{ \sqrt{n}  }   \ \ \    \ \ \ \ \ and \ \ \    \ \ \ \ \textbf{W}(  S_{n_0,n}) \ll  \displaystyle \frac{ \ln n  }{ \sqrt{n}  }.
	\end{eqnarray}
\end{itemize}
Moreover,    (\ref{rdffs1}) and  (\ref{rdffs2})  remain   valid when $ S_{n_0,n}$ is replaced  by $-S_{n_0,n}$.
\end{theorem}

\section{Proofs of the main results}\label{sec4}
\setcounter{equation}{0}

For the sake of simplicity,   in the sequel we denote
$\rho_{n,k },$ $\overline{\rho}_{n,k} $, $\upsilon_{n,k}$ and $\tau_{n,k}$  by $\rho_{k },$ $\overline{\rho}_{k} $, $\upsilon_{k}$ and $\tau_{k}$, respectively.

\subsection{Proof of Theorem \ref{th4.s1}}

The proof of Theorem \ref{th4.s1} is a refinement of Lemma 3.3 of Grama and Haeusler \cite{GH00}.  Compared to the proof of   Grama and Haeusler \cite{GH00}, our proof does not need the assumption $\| V_n^2/ s_n^2-1 \|_{\infty} \rightarrow 0.$
In the proof of Theorem \ref{th4.s1}, we make use of the following   technical lemma   of Bolthausen  \cite{2} (cf.\, Lemma 1 therein),
which plays an important role.
\begin{lemma}
\label{LEMMA-APX-1}Let $X$ and $Y$ be two random variables. Then
\[
  \mathbf{K}(X)
\leq c_1  \mathbf{K}(X+Y) + c_2 \big\| \mathbf{E}\left[ Y^2|X\right] \big\| _\infty ^{1/2}.
\]
\end{lemma}

Now we are in the position to prove Theorem \ref{th4.s1}.
For $a>0,$ define  $s_{n+1}^2 = s_n^2.$
  Thus, by the definition, we have $\tau_{1}=\sqrt{s_n^2+a^2}$ and $\tau_{n+1}= a $.
For  $u, x\in \mathbf{R}$  and $y > 0,$ denote
\begin{equation}
\Phi _u(x,y)=\Phi \bigg( \frac{u-x}{ \sqrt{y} }  \bigg) .  \label{RA-7}
\end{equation}
Let $\mathcal{N}$ be  independent of the martingale $\textit{\textbf{X}}$.
By Lemma \ref{LEMMA-APX-1}, for any $a>0,$
we can deduce that
\begin{eqnarray}
\mathbf{K}(S_n/s_n)
 &\leq& c_1 \, \mathbf{K}(S_n/s_n + a\mathcal{N}/s_n )
+c_2\, \frac{a}{s_n } \nonumber  \\
&=&  c_1\sup_{u \in \mathbf{R}}\Big|
\mathbf{E} [\Phi _u(S_n/s_n, a^2 /s_n^2)]- \Phi (u)\Big|
+c_2\, \frac{a}{s_n } \nonumber\\
&\leq& c_1\sup_{u \in \mathbf{R}}\Big| \mathbf{E} [\Phi _u(S_n/s_n, a^2 /s_n^2 )]-\mathbf{E} [\Phi _u(S_0/s_n , \tau^2_1/s_n^2)]\Big|  \nonumber \\
&& \ +\, c_1\sup_{u \in \mathbf{R}}\Big| \mathbf{E} [\Phi _u(S_0/s_n, \tau^2_1/s_n^2)]-\Phi (u)\Big| +c_2\,  \frac{a}{s_n } \nonumber\\
&=& c_1\sup_{u \in \mathbf{R}}\Big| \mathbf{E} [\Phi _u(S_n/s_n, \tau_{n+1}^2 /s_n^2)]-\mathbf{E} [\Phi _u(S_0/s_n,  \tau^2_1/s_n^2)]\Big|  \  \nonumber \\
 && +\, c_1\sup_{u \in \mathbf{R}}\bigg| \Phi \bigg(\frac{u}{\sqrt{a^2/s_n^2+ 1}}  \bigg)-\Phi (u)\bigg| +c_2\,  \frac{a}{s_n } . \label{dfafs}
\end{eqnarray}
Notice that $ a/s_n \rightarrow 0$ as  $n\rightarrow \infty.  $
It is easy to see that
\begin{eqnarray}
 \bigg| \Phi \bigg(\frac{u}{\sqrt{a^2/s_n^2+ 1}}  \bigg)-\Phi (u)\bigg| & \leq &  c_3 \bigg|\frac{1}{\sqrt{a^2/s_n^2+ 1}}  -1 \bigg| \nonumber \\
 &\leq &  c_4\frac{a}{s_n } .
\end{eqnarray}
Therefore, by the last inequality  and (\ref{dfafs}),  we can deduce that
\begin{eqnarray}
 \mathbf{K}(S_n/s_n)  \,  \leq \,  c_1\,\sup_{u \in \mathbf{R}}\bigg| \mathbf{E} [\Phi _u(S_n/s_n, \tau_{n+1}^2 /s_n^2)]-\mathbf{E} [\Phi _u(S_0/s_n,  \tau^2_1/s_n^2)]\bigg|   +c_5\, \frac{a}{s_n }  .   \label{dsdff}
\end{eqnarray}
In the sequel, we give an estimation for the term $  \mathbf{E} [\Phi _u(S_n/s_n, \tau_{n+1}^2 /s_n^2)]-\mathbf{E} [\Phi _u(S_0/s_n,  \tau^2_1/s_n^2)]$.
By a simple telescoping sum, we have
\begin{eqnarray*}
 && \mathbf{E} [\Phi _u(S_n/s_n, \tau^2_{n+1} /s_n^2)]-\mathbf{E} [\Phi _u(S_0/s_n,  \tau^2_1/s_n^2)] \\
 &&\ \ \ \ \ \ \ \ \ \ \ \ \ \   = \ \sum_{k=1}^n\mathbf{E}
\Big[ \Phi _u(S_k/s_n, \tau^2_{k+1} /s_n^2)-\Phi _u(S_{k-1}/s_n, \tau^2_{k} /s_n^2)  \Big] .
\end{eqnarray*}
Denote $$\xi_k= \frac{X_k}{s_n }   , \ \ \ \ \  \  1 \leq k \leq n.$$
Using the fact that
\[
\frac{\partial ^2}{\partial x^2}\Phi _u(x,y)=2\frac \partial {\partial
y}\Phi _u(x,y)
\]
and
that   $\frac{\partial ^2}{\partial x^2}\Phi
_u(S_{k-1}/s_n, \tau^2_{k+1} /s_n^2)$ is measurable with  respect to $\mathcal{F}_{k-1}$,
we obtain the following equality
\begin{equation}\nonumber
\mathbf{E} [\Phi _u(S_n/s_n, \tau^2_{n+1} /s_n^2)]-\mathbf{E} [\Phi _u(S_0/s_n,  \tau^2_1/s_n^2)]\ =\ I_1+I_2-I_3,
\end{equation}
where
\begin{eqnarray}
I_1\ &=&\  \sum_{k=1}^n \mathbf{E} \bigg[  \frac{}{} \Phi _u(S_k/s_n, \tau^2_{k+1} /s_n^2)-\Phi _u(S_{k-1}/s_n, \tau^2_{k+1} /s_n^2) \nonumber \\
&& \ \ \ \ \ \ \ \ \  \    +\, \frac \partial {\partial x}\Phi
_u(S_{k-1}/s_n, \tau^2_{k+1} /s_n^2)\xi_k-\frac 12\frac{\partial ^2}{\partial x^2}\Phi
_u(S_{k-1}/s_n, \tau^2_{k+1} /s_n^2)\xi _k^2   \bigg] ,   \nonumber\\
I_2 \ &=&\  \frac 12 \sum_{k=1}^n \mathbf{E}   \bigg[\frac{\partial ^2}{\partial x^2}\Phi
_u(S_{k-1}/s_n, \tau^2_{k+1} /s_n^2)\Big( \mathbf{E} [ \xi _k^2   | \mathcal{F}_{k-1}]   - \overline{\sigma}_k^2 /s_n^2  \Big) \bigg] ,\quad \quad \
 \nonumber \\
I_3\ &=&\ \sum_{k=1}^n\mathbf{E}   \bigg[ \Phi _u(S_{k-1}/s_n, \tau^2_{k} /s_n^2)-\Phi
_u(S_{k-1}/s_n, \tau^2_{k-1} /s_n^2)-\frac \partial {\partial y}\Phi _u(S_{k-1}/s_n, \tau^2_{k-1} /s_n^2)\frac{ \overline{\sigma}_k^2 }{s_n^2}
\bigg].  \nonumber
\end{eqnarray}
From (\ref{dsdff}) and the definitions of $I_1, I_2$ and $I_3$, we can deduce that
\begin{eqnarray}
 \mathbf{K}(S_n/s_n)  \leq  C\, \Big(  |I_1|+|I_2|+|I_3| + \frac{a}{s_n } \Big).
\end{eqnarray}
In the sequel,  we  give some estimates for the terms $|I_1|,$ $|I_2|$ and $|I_3|.$

\textbf{a)}  Control of $|I_1|.$
It is easy to see that
for all $1\leq k \leq n,$
\begin{eqnarray}
&&  R_k:=  \frac{}{} \Phi _u(S_k/s_n, \tau^2_{k+1} /s_n^2)-\Phi _u(S_{k-1}/s_n, \tau^2_{k+1} /s_n^2) \, +\, \frac \partial {\partial x}\Phi
_u(S_{k-1}/s_n, \tau^2_{k+1} /s_n^2)\, \xi_k   \nonumber \\
&& \ \ \ \ \ \ \ \ \  \ -\frac 12\frac{\partial ^2}{\partial x^2}\Phi
_u(S_{k-1}/s_n, \tau^2_{k+1} /s_n^2)\, \xi _k^2    \nonumber  \\
&& \ \ \ \ \,    =\   \Phi \Big( H_{k-1} -  \frac{ X_k}{\tau_{k+1}} \Big)-\Phi( H_{k-1})
    +\, \Phi'(H_{k-1}) \frac{ X_k}{\tau_{k+1}} -\frac 12\Phi''(H_{k-1}) \Big(\frac{ X_k}{\tau_{k+1}}\Big)^2 , \nonumber
\end{eqnarray}
where
$$\displaystyle H_{k-1}= \frac{u-S_{k-1}/s_n}{\tau_{k+1} /s_n} .$$
 When   $|X_k|\leq \tau_{k+1}$, it is easy to see that
\begin{eqnarray}
 \Big| \mathbf{E}[R_k \mathbf{1}_{ \{|X_k|\leq \tau_{k+1}  \}} ]   \Big|      &\leq &  \bigg|  \mathbf{E}\Big[  \frac1{6} \Phi^{(3)}  \Big( H_{k-1} - \vartheta \frac{ X_k}{\tau_{k+1}} \Big) \Big(\frac{ X_k}{\tau_{k+1}}\Big)^3 \textbf{1}_{\{|X_k|\leq \tau_{k+1}\}}   \Big] \bigg|\nonumber  \\
&\leq &   \frac{C  }{\tau_{k+1}^3}  \mathbf{E}\big[    |X_k|^3  \textbf{1}_{\{|X_k|\leq \tau_{k+1}\}} \big], \label{dfs01}
\end{eqnarray}
where  $\vartheta$ stands for some values or
random variables satisfying $0 \leq \vartheta \leq 1$, which may represent different values at different places.
Next, we consider the case  $|X_k|> \tau_{k+1}$.  It is easy to see that  for all $|\Delta x|>1 ,$
\begin{eqnarray*}
 \bigg|\Phi(x+\Delta x)-\Phi(x) - \Phi'(x)  \Delta x - \frac12 \Phi''(x) (\Delta x)^2 \bigg|
\ \leq \  \frac{c_2}{ 2+  x^2  }   \, |\Delta x|^{2 }.
\end{eqnarray*}
Therefore, we have
\begin{eqnarray}\nonumber
 \Big|  R_k \mathbf{1}_{ \{|X_k|> \tau_{k+1}  \}}    \Big|   \leq   G(H_{k-1}) \, \Big(\frac{ X_k}{\tau_{k+1}}\Big)^2\mathbf{1}_{ \{|X_k|> \tau_{k+1}  \}}  ,
\end{eqnarray}
where
$ \displaystyle G(z)=\frac{c_2}{ 2+  z^{2}  } .$
As $G$ is bounded,  we deduce that
\begin{eqnarray}\label{dfs02}
 \Big|\mathbf{E}[R_k \mathbf{1}_{ \{|X_k|> \tau_{k+1}  \}} ]  \Big| \leq     \frac{C  }{\tau_{k+1}^2}  \mathbf{E}\big[    X_k ^2  \textbf{1}_{\{|X_k|> \tau_{k+1}\}} \big].
\end{eqnarray}

Applying (\ref{dfs01}) and (\ref{dfs02}) to the estimation of $|I_1|,$ we get
\begin{eqnarray}
|I_1| &\leq& \sum_{k=1}^n \Big|\mathbf{E}R_k \Big|  \leq \sum_{k=1}^n \Big|\mathbf{E}[R_k \mathbf{1}_{ \{|X_k|\leq \tau_{k+1}  \}} ]   \Big| +  \sum_{k=1}^n \Big| \mathbf{E}[R_k \mathbf{1}_{ \{|X_k|>\tau_{k+1}  \}} ]   \Big|  \nonumber \\
 &\leq&  C  \,\sum_{k=1}^n \Bigg(  \frac{ \mathbf{E}\big[  |X_k|^3  \textbf{1}_{\{|X_k|\leq \tau_{k+1}\}} \big]   }{\tau_{k+1}^3}+  \frac{ \mathbf{E}\big[  X_k ^2  \textbf{1}_{\{|X_k|> \tau_{k+1}\}} \big]  }{\tau_{k+1}^2}  \Bigg) .
\end{eqnarray}

\textbf{b)}  Control of $|I_2|.$ For $|I_2|,$ we have
\begin{eqnarray*}
\left|
I_2\right| &\leq& \Bigg|\frac{1}{\tau_{k+1}^2 /s_n^2}\sum_{k=1}^n  \mathbf{E}   \bigg[  \varphi'(H_{k-1}) \Big(\mathbf{E} [ \xi _k^2   | \mathcal{F}_{k-1}]   - \overline{\sigma}_k^2 /s_n^2  \Big) \bigg] \Bigg|\\
&\leq& C\,\sum_{k=1}^n  \frac{1}{\tau_{k+1}^2 } \mathbf{E}   \big |  \sigma_k^2  - \overline{\sigma}_k^2       \big|,
\end{eqnarray*}
where  $\varphi$ is  the density function of the standard normal random variable.

 \textbf{c)}   Control of $|I_3|.$ For $|I_3|,$ by a two-term Taylor's expansion, it follows that
\[
I_3=\frac{1}{8}\,\sum_{k=1}^n\mathbf{E}   \Bigg[\frac 1{\tau_{k }^2/s_n^2-\vartheta   \overline{\sigma}^2_k/s_n^2  }  \varphi
^{\prime \prime \prime }\bigg( \frac{u-X_{k-1}/s_n}{\sqrt{ \tau_{k }^2/s_n^2-\vartheta   \overline{\sigma}^2_k/s_n^2}\   }
\bigg)  \Big(\frac{\overline{\sigma}^2_k }{s_n^2 }  \Big)^2\Bigg],
\]
where  $0 \leq \vartheta \leq 1$.
As $ \varphi^{\prime \prime \prime }$ is bounded,    we obtain
\begin{eqnarray}\label{fgfdd}
\left|I_3\right|  \ \leq \  C \,  \sum_{k=1}^n\frac 1{  \tau_{k-1}^2 /s_n^2 } \ \Big(\frac{\overline{\sigma}^2_k }{s_n^2 }  \Big)^2    .
\end{eqnarray}

Applying the upper bounds of $|I_1|, |I_2|$ and $|I_3|$ to (\ref{dsdff}),  we get
\begin{eqnarray}\label{ds5sfs}
 \mathbf{K}(S_n/s_n)
   \!  &\leq& \! C  \bigg[  \,\sum_{k=1}^n \bigg(  \frac{ \mathbf{E}\big[  |X_k|^3  \textbf{1}_{\{|X_k|\leq \tau_{k+1}\}} \big]   }{\tau_{k+1}^3}+  \frac{ \mathbf{E}\big[  X_k ^2  \textbf{1}_{\{|X_k|> \tau_{k+1}\}} \big] \ +\  \mathbf{E}   \big |  \sigma_k^2  - \overline{\sigma}_k^2       \big| }{\tau_{k+1}^2}   \nonumber  \\
  && \ \ \ \ \ \ \    + \ \frac 1{  \tau_{k-1}^2 /s_n^2 } \ \Big(\frac{\overline{\sigma}^2_k }{s_n^2 }  \Big)^2 \bigg)   +      \frac{a}{s_n } \ \bigg].
\end{eqnarray}
Notice that $\overline{\sigma}^2_k / s_n^2   \leq  || \textit{\textbf{X}} ||_2^2/ s_n^2  \rightarrow 0$ as $n \rightarrow \infty.$
Thus, when $a\geq1,$ we have
\begin{eqnarray}
 \sum_{k=1}^n \frac 1{  \tau_{k-1}^2 /s_n^2 } \ \Big(\frac{\overline{\sigma}^2_k }{s_n^2 }  \Big)^2 \  &\leq & \ \sum_{k=1}^n \frac 1{  1 - s_{k-1}^2/ s_n^2  +a^2/ s_n^2  } \  \frac{\overline{\sigma}^2_k }{s_n^2 }   \ \frac{|| \textit{\textbf{X}} ||_2^2}{s_n^2 }  \nonumber \\
 &\leq & \ C  \, \frac{|| \textit{\textbf{X}} ||_2^2}{s_n^2 } \ln \Big(2+  \frac{s_n^2}{a^2} \Big) \nonumber \\
  &\leq & \ C  \,   \frac{a}{s_n }. \label{fsdh5}
\end{eqnarray}
Applying the last inequality to  (\ref{ds5sfs}), we get the  desired inequality.
\hfill\qed

\subsection{Proof of Corollary   \ref{co2.5} }
Notice that $\overline{\sigma}^2_k / s_n^2   \leq  || \textit{\textbf{X}} ||_2^2/ s_n^2  \rightarrow 0$ as $n \rightarrow \infty.$
Thus, we have
\begin{eqnarray}
\sum_{k=1}^n \frac{ \mathbf{E}\big[  |X_k|^3  \textbf{1}_{\{|X_k|\leq \tau_{k+1}\}} \big]   }{\tau_{k+1}^3} &\leq& \sum_{k=1}^n\frac{ \mathbf{E}  |X_k|^{2+ \delta}  }{ \tau_{k+1}^{2+\delta} }  \ \leq\ C\, \sum_{k=1}^n\frac{  \overline{\sigma}_k^2 }{ \tau_{k+1}^{2+\delta} } \nonumber\\
 &=& \  C\,\sum_{k=1}^n\frac{ 1}{ (\tau_{k+1}/s_n)^{2+\delta} }  \frac{\overline{\sigma}_k^2}{s_n^2} \frac{1}{s_n^{\delta}}\nonumber\\
&\leq& \  C\, \frac{1}{a^\delta }  .
\end{eqnarray}
Similarly, we have
\begin{eqnarray}
\sum_{k=1}^n \frac{ \mathbf{E}\big[  |X_k|^2  \textbf{1}_{\{|X_k|> \tau_{k+1}\}} \big]   }{\tau_{k+1}^2} &\leq& \sum_{k=1}^n\frac{ \mathbf{E}  |X_k|^{2+ \delta}  }{ \tau_{k+1}^{2+\delta} }  \ \leq\ C\,   \frac{1}{a^\delta }  .
\end{eqnarray}
Using the condition
$\mathbf{E} \big|\sigma_k^2-\overline{\sigma}_k^2 \big|\leq C \overline{\sigma}_k^2 s_n^{-\alpha},$ we get
\begin{eqnarray}\label{une4.3}
\sum_{k=1}^n\frac{ \mathbf{E}   \big |  \sigma_k^2  - \overline{\sigma}_k^2       \big| }{\tau_{k+1}^2}&\leq& C \sum_{k=1}^n\frac{1}{ \tau_{k+1}^{2 } }  \frac{\overline{\sigma}_k^2}{s_n^{\alpha}} \ \leq\ C\,   \frac{\ln(2+ s_n^2/a^2)}{s_n^{\alpha}}  .
\end{eqnarray}
Thus, by Theorem  \ref{th4.s1}, we deduce that for all $a\geq 1,$
\begin{eqnarray*}
 \mathbf{K}(S_n/s_n)
   \   \leq  \  C  \bigg(   \frac{1}{a^\delta }  \ +\    \frac{\ln(2+ s_n^2/a^2)}{s_n^{\alpha}}\ + \ \frac{a}{s_n  }  \bigg)   .
\end{eqnarray*}
Taking $a= s_n^{1/(1+\delta)}$ in the last inequality, we get
\begin{eqnarray*}
 \mathbf{K}(S_n/s_n)
   \   \leq  \   C \, \bigg(\displaystyle \frac{ 1}{ s_n^{\delta /(1+ \delta) }}+ \frac{\ln s_n  }{ s_n^{\alpha } }  \bigg)  ,
\end{eqnarray*}
which gives the desired inequality.

\subsection{Proof of Corollary   \ref{th4.21}  }
By the conditions of Corollary   \ref{th4.21}, it is easy to see that  $\sigma^2 \leq  \bar{\sigma}_k^2  \leq ||\mathbf{X}||^2_{2+\delta}.$
By the fact $\mathbf{E} \big|\sigma_k^2-\overline{\sigma}_k^2 \big|=O( k^{-\alpha}) \   (k\rightarrow \infty)$, we can deduce that
for $\alpha \in (0, 1) $ and all $n\geq 2,$
  \begin{eqnarray*}
\sum_{k=1}^n\frac{ \mathbf{E}   \big |  \sigma_k^2  - \overline{\sigma}_k^2       \big| }{\tau_{k+1}^2}&\leq&  C\, \sum_{k=1}^n  \frac{ 1 }{  n -k + 1   \ } k^{-\alpha} \\
  &=& C\, \sum_{k=1}^n  \frac{ 1 }{  1 -k/n + 1/n   \ } \frac{1}{ (k/n)^{\alpha}} \, \frac{1}{n} \, n^{ -\alpha} \\
  &\leq& C\, \int_{0}^1 \frac{}{}\frac{}{}   \frac{ 1 }{   1 -x +1/n\ } \frac{1}{ x^{\alpha}} \,dx \, n^{ -\alpha} \\
  &\leq& C\,  \frac{\ln  n  }{ n^{\alpha } }.
\end{eqnarray*}
The remain of the proof is similar to the proof of  Corollary   \ref{co2.5}.

\subsection{Proof of Theorem \ref{th4.2}}

For the proof of Theorem \ref{th4.2}, we make use of the following technical lemma of Bolthausen  (cf.\ Lemma 2  of \cite{2}).

\begin{lemma}
\label{LEMMA-APX-2}Let $G(x)$ be an integrable function on $\mathbf{\mathbf{R}}$ of bounded variation $||G||_V$,
$X$   a random variable and $a,$ $b\neq 0$   real numbers. Then
\[
\mathbf{E}\Big[ \, G\Big( \frac{X+a}b\Big) \Big] \leq ||G||_V \, \mathbf{K}(X) + ||G||_1 \, |b|,
\]
where  $||G||_1$ is the $L_1(\mathbf{R})$ norm of $G(x).$
\end{lemma}

It is easy to see that
for all $1\leq k \leq n,$
\begin{eqnarray}
&&  R_k    =\   \Phi \Big( H_{k-1} -  \frac{ X_k}{\tau_{k+1}} \Big)-\Phi( H_{k-1})
    +\, \Phi'(H_{k-1}) \frac{ X_k}{\tau_{k+1}} -\frac 12\Phi''(H_{k-1}) \Big(\frac{ X_k}{\tau_{k+1}}\Big)^2 . \nonumber
\end{eqnarray}
For the estimation of $ |\mathbf{E} R_k|,$
we distinguish two cases as follows.

\textbf{a)}  Case $|X_k|\leq \tau_{k+1}$.
When $|X_k|\leq \tau_{k+1}$, it is easy to see that
\begin{eqnarray}
 \Big| \mathbf{E}[R_k \mathbf{1}_{ \{|X_k|\leq \tau_{k+1}  \}} ]   \Big|      &\leq &  \bigg|  \mathbf{E}\Big[  \frac1{6} \Phi^{(3)}  \Big( H_{k-1} - \vartheta \frac{ X_k}{\tau_{k+1}} \Big) \Big(\frac{ X_k}{\tau_{k+1}}\Big)^3 \textbf{1}_{\{|X_k|\leq \tau_{k+1}\}}   \Big] \bigg|\nonumber  \\
&\leq &   \frac{C  }{\tau_{k+1}^3}  \mathbf{E}\big[ H(H_{k-1})   |X_k|^3  \textbf{1}_{\{|X_k|\leq \tau_{k+1}\}} \big]  \nonumber  \\
&\leq&  \frac{C  }{\tau_{k+1}^3} \big\| \mathbf{E} [  |X_k|^3  \mathbf{1}_{\{|X_k|\leq \tau_{k+1}\}} |\mathcal{F}_{k-1} ]  \big\|_{\infty} \mathbf{E}[  H(H_{k-1})],   \label{dfssf01}
\end{eqnarray}
where $H(x)=\sup_{|t|\leq 1}   |\Phi^{(3)}(x+t) |$ and $\vartheta$ stands for some values or
random variables satisfying $0 \leq \vartheta \leq 1$, which may represent different values at different places.
By Lemma \ref{LEMMA-APX-2}, we have
\begin{eqnarray}
\mathbf{E}[  H(H_{k-1})]  \leq    C \, \mathbf{K}(S_{k-1}/s_n) + C \, |\tau_{k+1} /s_n  |.
\end{eqnarray}
By Lemma \ref{LEMMA-APX-1}, we deduce that for $a$ large enough,
\begin{eqnarray}
 \mathbf{E}[  H(H_{k-1})]
& \leq &
 C \, \mathbf{K}(S_{n}/s_n) + C \big\| \rho_{k}  / s_n   \big\| _\infty + C \, |\tau_{k+1} /s_n  | \nonumber \\
& \leq & C \, \mathbf{K}(S_{n}/s_n) +  C \, |\tau_{k+1} /s_n  | ,   \nonumber
\end{eqnarray}
where the last line using the fact that $\rho_{n,k}^2 \leq  C\, \overline{\rho}^2_{n,k} $.  Applying the last line to (\ref{dfssf01}), we get
\begin{eqnarray}
 \Big| \mathbf{E}[R_k \mathbf{1}_{ \{|X_k|\leq \tau_{k+1}  \}} ]   \Big|      &\leq &    \frac{C \big\| \mathbf{E} [  |X_k|^3  \mathbf{1}_{\{|X_k|\leq \tau_{k+1}\}} |\mathcal{F}_{k-1} ]  \big\|_{\infty} }{\tau_{k+1}^3}  \mathbf{K}(S_{n}/s_n)  \nonumber \\
  && \ \ \ \ \ \ \ \ \ \ \ + \ \frac{C  }{s_n\,} \frac{\big\| \mathbf{E} [  |X_k|^3  \mathbf{1}_{\{|X_k|\leq \tau_{k+1}\}} |\mathcal{F}_{k-1} ]  \big\|_{\infty}}{\tau_{k+1}^2}
. \label{dfssdd1}
\end{eqnarray}

\textbf{b)}  Case $|X_k|>\tau_{k+1}$.    It is easy to see that  for all $|\Delta x|>1 ,$
\begin{eqnarray*}
 \bigg|\Phi(x+\Delta x)-\Phi(x) - \Phi'(x)  \Delta x - \frac12 \Phi''(x) (\Delta x)^2 \bigg|
\ \leq \  \frac{c_2}{ 2+  x^2  }   \, |\Delta x|^{2 }.
\end{eqnarray*}
Therefore, we have
\begin{eqnarray}\nonumber
 \Big| \mathbf{E}[R_k \mathbf{1}_{ \{|X_k|> \tau_{k+1}  \}}  ]   \Big|   \leq   G(H_{k-1}) \, \Big(\frac{ X_k}{\tau_{k+1}}\Big)^2\mathbf{1}_{ \{|X_k|> \tau_{k+1}  \}}  ,
\end{eqnarray}
where
$ \displaystyle G(z)=\frac{c_2}{ 2+  z^{2}  } .$
As $H_{k-1}$ is $\mathcal{F}_{k-1}$-measurable,  we deduce that
\begin{eqnarray}
 \Big|\mathbf{E}[R_k \mathbf{1}_{ \{|X_k|> \tau_{k+1}  \}} ]  \Big| &\leq&     \frac{C  }{\tau_{k+1}^2}  \mathbf{E}\big[ G(H_{k-1})  X_k ^2  \textbf{1}_{\{|X_k|> \tau_{k+1}\}} \big] \nonumber \\
 &\leq&  \frac{C  }{\tau_{k+1}^2} \big\| \mathbf{E} [  X_k^2  \mathbf{1}_{\{|X_k|> \tau_{k+1}\}} |\mathcal{F}_{k-1} ]  \big\|_{\infty} \mathbf{E}[  G(H_{k-1})].\label{dfgfs02}
\end{eqnarray}
By an argument similar to the proof of (\ref{dfssdd1}), we deduce that
\begin{eqnarray}
 \Big| \mathbf{E}[R_k \mathbf{1}_{ \{|X_k|> \tau_{k+1}  \}} ]   \Big|      &\leq &    \frac{C  \big\| \mathbf{E} [  X_k^2  \mathbf{1}_{\{|X_k|> \tau_{k+1}\}} |\mathcal{F}_{k-1} ]  \big\|_{\infty}}{\tau_{k+1}^2}  \mathbf{K}(S_{n}/s_n)  \nonumber \\
  &&\ \ \ \ \ \ \ \ \ \ \ +\ \frac{C  }{s_n\,} \frac{\big\| \mathbf{E} [   X_k^2  \mathbf{1}_{\{|X_k|> \tau_{k+1}\}} |\mathcal{F}_{k-1} ]  \big\|_{\infty}}{\tau_{k+1} } \label{dfssdd2}
.
\end{eqnarray}

Applying (\ref{dfssdd1}) and (\ref{dfssdd2}) to the estimation of $|I_1|,$ we get
\begin{eqnarray}
|I_1| &\leq &\sum_{k=1}^n\big|\mathbf{E} R_k    \big|  \leq \sum_{k=1}^n \Big| \mathbf{E}[R_k \mathbf{1}_{ \{|X_k|\leq \tau_{k+1}  \}} ]   \Big|   +\sum_{k=1}^n\Big| \mathbf{E}[R_k \mathbf{1}_{ \{|X_k|> \tau_{k+1}  \}} ]   \Big|  \nonumber \\
&\leq&  C \sum_{k=1}^n\bigg[  \frac{ \big\| \mathbf{E} [  |X_k|^3  \mathbf{1}_{\{|X_k|\leq \tau_{k+1}\}} |\mathcal{F}_{k-1} ]  \big\|_{\infty}}{\tau_{k+1}^3}  +  \frac{ \big\| \mathbf{E} [  X_k^2  \mathbf{1}_{\{|X_k|> \tau_{k+1}\}} |\mathcal{F}_{k-1} ]  \big\|_{\infty}}{\tau_{k+1}^2}  \bigg] \mathbf{K}(S_{n}/s_n)  \nonumber \\
&& \ +\ \frac{C  }{s_n\,} \sum_{k=1}^n\bigg[  \frac{ \big\| \mathbf{E} [  |X_k|^3  \mathbf{1}_{\{|X_k|\leq \tau_{k+1}\}} |\mathcal{F}_{k-1} ]  \big\|_{\infty}}{\tau_{k+1}^2}  +  \frac{ \big\| \mathbf{E} [  X_k^2  \mathbf{1}_{\{|X_k|> \tau_{k+1}\}} |\mathcal{F}_{k-1} ]  \big\|_{\infty}}{\tau_{k+1} }  \bigg].
\end{eqnarray}
By the controls of $|I_2|$ and $|I_3|$ in the proof of Theorem \ref{th4.s1}, we deduce that
\begin{eqnarray*}
 \mathbf{K}(S_n/s_n)
     &\leq&  C_* \sum_{k=1}^n\bigg[  \frac{ \big\| \mathbf{E} [  |X_k|^3  \mathbf{1}_{\{|X_k|\leq \tau_{k+1}\}} |\mathcal{F}_{k-1} ]  \big\|_{\infty}}{\tau_{k+1}^3}  +  \frac{ \big\| \mathbf{E} [  X_k^2  \mathbf{1}_{\{|X_k|> \tau_{k+1}\}} |\mathcal{F}_{k-1} ]  \big\|_{\infty}}{\tau_{k+1}^2}  \bigg] \mathbf{K}(S_{n}/s_n)  \nonumber \\
&& \ +\ \frac{C  }{s_n\,} \sum_{k=1}^n\bigg[  \frac{ \big\| \mathbf{E} [  |X_k|^3  \mathbf{1}_{\{|X_k|\leq \tau_{k+1}\}} |\mathcal{F}_{k-1} ]  \big\|_{\infty}}{\tau_{k+1}^2}  +  \frac{ \big\| \mathbf{E} [  X_k^2  \mathbf{1}_{\{|X_k|> \tau_{k+1}\}} |\mathcal{F}_{k-1} ]  \big\|_{\infty}}{\tau_{k+1} }  \bigg]\\
&& +\  C\, \Bigg[ \sum_{k=1}^n \frac{ \mathbf{E}   \big |  \sigma_k^2  - \overline{\sigma}_k^2       \big| }{\tau_{k+1}^2} \ +  \ \frac 1{  \tau_{k-1}^2 /s_n^2 } \ \Big(\frac{\overline{\sigma}^2_k }{s_n^2 }  \Big)^2    \ \Bigg]  +      \frac{a}{s_n }.
\end{eqnarray*}
By condition (\ref{fsdfzn}), for all $a$ large enough, we have
$$\sum_{k=1}^n\bigg[  \frac{ \big\| \mathbf{E} [  |X_k|^3  \mathbf{1}_{\{|X_k|\leq \tau_{k+1}\}} |\mathcal{F}_{k-1} ]  \big\|_{\infty}}{\tau_{k+1}^3}  +  \frac{ \big\| \mathbf{E} [  X_k^2  \mathbf{1}_{\{|X_k|> \tau_{k+1}\}} |\mathcal{F}_{k-1} ]  \big\|_{\infty}}{\tau_{k+1}^2}  \bigg] \leq \frac{1}{2C_*},$$
and,  by (\ref{fsdh5}),
$$\frac 1{  \tau_{k-1}^2 /s_n^2 } \ \Big(\frac{\overline{\sigma}^2_k }{s_n^2 }  \Big)^2 \leq  C     \frac{a}{s_n }. $$
Therefore, it holds for all $a$ large enough,
\begin{eqnarray*}
 \mathbf{K}(S_n/s_n)
     &\leq& \frac{C  }{s_n\,} \sum_{k=1}^n\bigg[  \frac{ \big\| \mathbf{E} [  |X_k|^3  \mathbf{1}_{\{|X_k|\leq \tau_{k+1}\}} |\mathcal{F}_{k-1} ]  \big\|_{\infty}}{\tau_{k+1}^2}  +  \frac{ \big\| \mathbf{E} [  X_k^2  \mathbf{1}_{\{|X_k|> \tau_{k+1}\}} |\mathcal{F}_{k-1} ]  \big\|_{\infty}}{\tau_{k+1} }  \bigg]\\
&& + \ \frac12 \mathbf{K}(S_n/s_n) \ +\  C\,  \bigg[  \sum_{k=1}^n \frac{ \mathbf{E}   \big |  \sigma_k^2  - \overline{\sigma}_k^2       \big| }{\tau_{k+1}^2}      +      \frac{a}{s_n } \bigg] ,
\end{eqnarray*}
which implies
the desired inequality.

\subsection{Proof of Corollary   \ref{fssdvs}   }
We only need to consider the case $\delta \in (0, 1].$
We first verify the condition (\ref{fsdfzn}).
Notice that $\overline{\sigma}^2_k / s_n^2   \leq  ||\mathbf{X}||_2/ s_n^2  \rightarrow 0$ as $n \rightarrow \infty.$
Thus,   we have
\begin{eqnarray}
\sum_{k=1}^n \frac{ \big\| \mathbf{E} [  |X_k|^3  \mathbf{1}_{\{|X_k|\leq \tau_{k+1}\}} |\mathcal{F}_{k-1} ]  \big\|_{\infty}}{\tau_{k+1}^3} &\leq& \sum_{k=1}^n\frac{\big\| \mathbf{E} [     |X_k|^{2+ \delta} |\mathcal{F}_{k-1} ] \big\|_{\infty} }{ \tau_{k+1}^{2+\delta} }  \ \leq\ C\, \sum_{k=1}^n\frac{  \overline{\sigma}_k^2 }{ \tau_{k+1}^{2+\delta} } \nonumber\\
 &\leq& C\, \sum_{k=1}^n\frac{  \overline{\sigma}_k^2 }{ (\tau_{k+1}/s_n)^{2+\delta/2} } \frac{  \overline{\sigma}_k^2 }{ s_n^2 } \frac{1}{s_n^{\delta/2}} \frac{ 1}{   a^{\delta/2} }\  \nonumber \\
  &\leq& C\, (s_n/a)^{\delta/2} \frac{1}{s_n^{\delta/2}} \frac{ 1}{   a^{\delta/2} }\ \nonumber \\
   &=& C\,   \frac{1 }{a^\delta \,}  \rightarrow 0.
\end{eqnarray}
Similarly, we have
\begin{eqnarray}
\sum_{k=1}^n \frac{ \big\| \mathbf{E} [  X_k^2  \mathbf{1}_{\{|X_k|> \tau_{k+1}\}} |\mathcal{F}_{k-1} ]  \big\|_{\infty}  }{\tau_{k+1}^2} &\leq& \sum_{k=1}^n\frac{ \mathbf{E}  |X_k|^{2+ \delta}  }{ \tau_{k+1}^{2+\delta} }  \ \leq\ C\, \frac{1 }{a^\delta\, }  \rightarrow 0   .
\end{eqnarray}
Thus the condition (\ref{fsdfzn}) is satisfied.
Taking $a=s_n,$   we have for $\delta \in (0, 1)$,
\begin{eqnarray}
\sum_{k=1}^n \frac{ \big\| \mathbf{E} [  |X_k|^3  \mathbf{1}_{\{|X_k|\leq \tau_{k+1}\}} |\mathcal{F}_{k-1} ]  \big\|_{\infty}}{\tau_{k+1}^2} &\leq& \sum_{k=1}^n\frac{\big\| \mathbf{E} [     |X_k|^{2+ \delta} |\mathcal{F}_{k-1} ] \big\|_{\infty} }{ \tau_{k+1}^{1+\delta} }  \ \leq\ C\, \sum_{k=1}^n\frac{  \overline{\sigma}_k^2 }{ \tau_{k+1}^{1+\delta} } \nonumber\\
 &=& \  C\,\sum_{k=1}^n\frac{ 1}{ (\tau_{k+1}/s_n)^{1+\delta} }  \frac{\overline{\sigma}_k^2}{s_n^2}  s_n^{1-\delta} \nonumber\\
&\leq& \  C\, s_n^{1-\delta}   . \label{fdsds1}
\end{eqnarray}
Similarly, we have for $\delta \in (0, 1)$
\begin{eqnarray}
\sum_{k=1}^n \frac{ \big\| \mathbf{E} [  X_k^2  \mathbf{1}_{\{|X_k|> \tau_{k+1}\}} |\mathcal{F}_{k-1} ]  \big\|_{\infty}  }{\tau_{k+1} } &\leq& \sum_{k=1}^n\frac{\big\| \mathbf{E} [     |X_k|^{2+ \delta} |\mathcal{F}_{k-1} ] \big\|_{\infty} }{ \tau_{k+1}^{1+\delta} } \ \leq\  C\, s_n^{1-\delta}  .\label{fdsds2}
\end{eqnarray}
Applying the  inequalities   (\ref{fdsds1}) and  (\ref{fdsds2}) to Theorem \ref{th4.2}, we obtain the desired inequality for  $\delta \in (0, 1)$.
For $\delta \geq 1,$ the inequality follows by a similar argument.

\subsection{Proof of Corollary   \ref{fsdvs}  }
The proof of  Corollary   \ref{fsdvs} is similar to the proof of Corollary  \ref{th4.21}.

For the sake of simplicity,   in the sequel we  denote
$\rho_{n,k },$ $\overline{\rho}_{n,k} $, $\upsilon_{n,k}$ and $\tau_{n,k}$  by $\rho_{k },$ $\overline{\rho}_{k} $, $\upsilon_{k}$ and $\tau_{k}$, respectively.
\subsection{Proof of Theorem \ref{theorem2.1}}
The proof is a modification of R\"ollin's argument.
Let $Z', Z_1,..., Z_n$ be a sequence of independent standard normal random variables, and they are independent of $\mathbf{X}$. Define
\begin{equation}
Z:=\sum_{i=1}^n \sigma_i Z_i, \qquad T_{k}:=\sum_{i=k}^n \sigma_i Z_i, \qquad k = 1, 2, ..., n.
\end{equation}
By the definition of  $Z$, $Z$ is a normal random variable with mean $0$ and variance $s_n^2$. As $\rho_k^2=s_n^2- V_k^2$ and $V_k^2$
is measurable with respect to $\mathcal{F}_{k-1}$,   it holds
\begin{equation}
\mathcal{L}(T_k| \mathcal{F}_{k-1}) \sim N(0,\rho_k^2),\
\end{equation}
where $\mathcal{L}(\cdot | \mathcal{F}_{k-1})$ represents the conditional distribution on $\mathcal{F}_{k-1}$.
Let $h$ be  a fixed 1-Lipschitz continuous  function. By the properties of  Lipschitz continuous  functions,
we have $||h'||_{\infty} \leq 1$, where $||\cdot||_{\infty}$ denotes the supremum norm with respect to Lebesgue's measure.
By the triangle inequality, it is easy to see that
	\begin{equation}
		\Big|\mathbf{E}[h(S_n)-h(Z)] \Big| \leq \Big |\mathbf{E}[h(S_n+aZ')-h(Z+aZ')] \Big|+2a.
	\end{equation}
Define $S_0=T_{n+1}=0.$
Using telescoping sum representation, we have
	\begin{equation}
		\mathbf{E}[h(S_n+aZ')-h(Z+aZ')] = \mathbf{E} \sum_{k=1}^n \mathbf{E}[R_k|\mathcal{F}_{k-1}],
	\end{equation}
where
	\begin{equation}
		R_k = h(S_k+T_{k+1}+aZ') - h(S_{k-1}+T_k+aZ').
	\end{equation}
Let $g$ be the unique bounded solution of the following equation
	\begin{equation}
		g'(x)-xg(x)=f(x)-\mathbf{E} f(Y),\ \ x \in \mathbf{R},
	\end{equation}
where $f(x) = h(tx+s)/t$, $s \in \mathbf{R}$, $t > 0$, $\mathcal{L}(Y) \sim N(0,1)$. Clearly, $f$ is a 1-Lipschitz continuous function.
Chen,  Goldstein and Shao~\cite{7} proved that
	\begin{equation}
		  ||g||_{\infty} \leq 2||f'||_{\infty}, \qquad ||g'||_{\infty} \leq \sqrt{\frac{2}{\pi}}||f'||_{\infty}, \qquad ||g''||_{\infty} \leq 2||f'||_{\infty}.
	\end{equation}
Set $f_{s,t}(w) := g((w - s)/t)$, where $w \in \mathbf{R}$. Then $f_{s,t}$ is also a Lipschitz function such that
	\begin{equation}\label{gfdgd}
		t^2 f'_{s,t}(w)-(w-s)f_{s,t}(w)=h(w)-\mathbf{E}h(tY+s), \quad w \in \mathbf{R}.
	\end{equation}
As $g$ is bounded, we have
	\begin{equation}\label{gsgs}
		||f_{s,t}||_{\infty} \leq 2||h'||_{\infty}, \qquad ||f'_{s,t}||_{\infty} \leq \frac{||h'||_{\infty}}{t}, \qquad ||f''_{s,t}||_{\infty} \leq \frac{2||h'||_{\infty}}{t^2}.
	\end{equation}
By the definition of $f_{s,t}$,  $f_{s,t}$, $f'_{s,t}$ and $f''_{s,t}$ can be regarded as measurable functions of $(s,t, w)$ from $\mathbf{R} \times \mathbf{R}^{+} \times \mathbf{R}$ to $\mathbf{R}$. Therefore,   $f_{U,V}(W)$, $f'_{U,V}(W)$ and $f''_{U,V}(W)$ can be regarded as random functions for the random variables $ U$, $V$ and  $W$, where $V>0$. Denote $T'_k:=T_k+aZ'$, then
\begin{equation}
	\mathcal{L}(T'_k|\mathcal{F}_{k-1})\sim N(0,\rho_k^2+a^2).
\end{equation}
Notice that $S_{k-1}$ and $ \sqrt{\rho_k^2+a^2}$ are $\mathcal{F}_{k-1}$-measurable. Thus
\begin{eqnarray*}	
\mathbf{E}[R_k|\mathcal{F}_{k-1}]&=&\mathbf{E}[h(S_k+T_{k+1}+aZ')-h(S_{k-1}+T_k+aZ')|\mathcal{F}_{k-1}]\\
&=&\mathbf{E}[h(S_k+T'_{k+1})-h(S_{k-1}+T'_k)|\mathcal{F}_{k-1}]\\
&=&\mathbf{E}[h(S_k+T'_{k+1})-h(S_{k-1}+v_k Y)|\mathcal{F}_{k-1}]\\
&=& \mathbf{E}[v_k^{2} f'_{S_{k-1},v_k}(S_k+T'_{k+1})-(S_k+T'_{k+1}-S_{k-1}) f_{S_{k-1},v_k}(S_k+T'_{k+1})|\mathcal{F}_{k-1}].
\end{eqnarray*}
The last line follows by  (\ref{gfdgd}). Since $v_k=\sqrt{\rho_k^2+a^2}$ and $\rho_k^2=\rho_{k+1}^2+\sigma_k^2$,	we have
	\begin{eqnarray}
		\mathbf{E}[R_k|\mathcal{F}_{k-1}]&=& \mathbf{E}[(\rho_k^2+a^2) f'_{S_{k-1},v_k}(S_k+T'_{k+1})-(X_k+T'_{k+1}) f_{S_{k-1},v_k}(S_k+T'_{k+1})|\mathcal{F}_{k-1}] \nonumber \\
		&=&\,\mathbf{E}[\sigma_k^2 f'_{S_{k-1},v_k}(S_k+T'_{k+1})-X_k f_{S_{k-1},v_k}(S_k+T'_{k+1})|\mathcal{F}_{k-1}] \label{rsffg}\\
		&&\    +\,\mathbf{E}[(\rho_{k+1}^2+a^2) f'_{S_{k-1},v_k}(S_k+T'_{k+1}) - T'_{k+1} f_{S_{k-1},v_k}(S_k+T'_{k+1})|\mathcal{F}_{k-1}]. \nonumber
	\end{eqnarray}
Notice that $f_{s,t}$~ is a Lipschitz function and $\mathcal{L}(T'_{k+1}|\mathcal{F}_{k-1}) \sim N(0,\rho_{k+1}^2+a^2)$.
For any normal random variable $Y$ with $\mathcal{L}(Y) \sim N(0,\sigma^2)$, if $\mathbf{E}{g'(Y)}$ exists, then the function $g$ satisfies the equation $\mathbf{E}[\sigma^2 g'(Y)-Yg(Y)]=0.$ Therefore, it holds
$$  \mathbf{E}[(\rho_{k+1}^2+a^2) f'_{S_{k-1},v_k}(S_k+T'_{k+1}) - T'_{k+1} f_{S_{k-1},v_k}(S_k+T'_{k+1})|\mathcal{F}_{k-1}]=0 .$$
Applying the last equality to (\ref{rsffg}), we get
\begin{eqnarray}
		\mathbf{E}[R_k|\mathcal{F}_{k-1}]
		 = \mathbf{E}[\sigma_k^2 f'_{S_{k-1},v_k}(S_k+T'_{k+1}) - X_k f_{S_{k-1},v_k}(S_k+T'_{k+1})|\mathcal{F}_{k-1}].
	\end{eqnarray}
Next, using  Taylor's expansions for $f'_{S_{k-1},v_k}(S_k+T'_{k+1})$ and $f_{S_{k-1},v_k}(S_k+T'_{k+1})$, we deduce that
	\begin{eqnarray}
\mathbf{E}[R_k|\mathcal{F}_{k-1}]
&=&  \mathbf{E}[ \sigma_k^2\mathbf{1}_{\{|X_k| < v_k\}} f'_{S_{k-1},v_k}(S_{k-1}+T'_{k+1})+
	\sigma_k^2 X_k\mathbf{1}_{\{|X_k| < v_k\}} f''_{S_{k-1},v_k}(S_{k-1}+\theta_1 X_k +T'_{k+1})  \nonumber\\
	&&    -X_k \mathbf{1}_{\{|X_k| < v_k\}}f_{S_{k-1},v_k}(S_{k-1}+T'_{k+1})-X_k^2\mathbf{1}_{\{|X_k| < v_k\}} f'_{S_{k-1},v_k}(S_{k-1} +T'_{k+1}) \nonumber \\
	&&    -X_k^3\mathbf{1}_{\{|X_k| < v_k\}} f''_{S_{k-1},v_k}(S_{k-1}+\theta_2 X_k +T'_{k+1})  \  \nonumber \\
&&    +  \sigma_k^2\mathbf{1}_{\{|X_k| \geq v_k\}} f'_{S_{k-1},v_k}(S_{k}+T'_{k+1})-X_k\mathbf{1}_{\{|X_k| \geq v_k\}} f_{S_{k-1},v_k}(S_{k-1}+T'_{k+1})  \nonumber  \\
	&&-X_k^2\mathbf{1}_{\{|X_k| \geq v_k\}} f'_{S_{k-1},v_k}(S_{k-1} +\theta_3 X_k+T'_{k+1}) |\mathcal{F}_{k-1}]      \nonumber \\
&=&  \mathbf{E}[ \sigma_k^2\mathbf{1}_{\{|X_k| < v_k\}} f'_{S_{k-1},v_k}(S_{k-1}+T'_{k+1})+
	\sigma_k^2 X_k\mathbf{1}_{\{|X_k| < v_k\}} f''_{S_{k-1},v_k}(S_{k-1}+\theta_1 X_k +T'_{k+1})  \nonumber \\
	&&    -X_k  f_{S_{k-1},v_k}(S_{k-1}+T'_{k+1})-X_k^2\mathbf{1}_{\{|X_k| < v_k\}} f'_{S_{k-1},v_k}(S_{k-1} +T'_{k+1}) \nonumber \\
	&&     -X_k^3\mathbf{1}_{\{|X_k| < v_k\}} f''_{S_{k-1},v_k}(S_{k-1}+\theta_2 X_k +T'_{k+1}) +  \sigma_k^2\mathbf{1}_{\{|X_k| \geq v_k\}} f'_{S_{k-1},v_k}(S_{k}+T'_{k+1}) \ \nonumber \\
&&    -X_k^2\mathbf{1}_{\{|X_k| \geq v_k\}} f'_{S_{k-1},v_k}(S_{k-1} +\theta_3 X_k+T'_{k+1})\, |\mathcal{F}_{k-1}]     \nonumber \\
&=& \mathbf{E}\Big[ \sigma_k^2  f'_{S_{k-1},v_k}(S_{k-1}+T'_{k+1})+
	\sigma_k^2 X_k\mathbf{1}_{\{|X_k| < v_k\}} f''_{S_{k-1},v_k}(S_{k-1}+\theta_1 X_k +T'_{k+1}) \nonumber \\
	&&   -X_k  f_{S_{k-1},v_k}(S_{k-1}+T'_{k+1})-X_k^2 f'_{S_{k-1},v_k}(S_{k-1} +T'_{k+1}) \nonumber \\
	&&    -X_k^3\mathbf{1}_{\{|X_k| < v_k\}} f''_{S_{k-1},v_k}(S_{k-1}+\theta_2 X_k +T'_{k+1})  \ \nonumber \\
&&    +  \sigma_k^2\mathbf{1}_{\{|X_k| \geq v_k\}} \Big(f'_{S_{k-1},v_k}(S_{k}+T'_{k+1}) - f'_{S_{k-1},v_k}(S_{k-1}+T'_{k+1}) \Big)  \nonumber \\
&&  -X_k^2\mathbf{1}_{\{|X_k| \geq v_k\}} \Big(f'_{S_{k-1},v_k}(S_{k-1} +\theta_3 X_k+T'_{k+1}) - f'_{S_{k-1},v_k}(S_{k-1}+T'_{k+1}) \Big) \   \Big |\mathcal{F}_{k-1} \Big] ,\nonumber
	\end{eqnarray}
where $0\leq\theta_1,\theta_2, \theta_3 \leq1$. By the last  equality and (\ref{gsgs}),  we deduce that
	\begin{eqnarray}
	\mathbf{E}[R_k|\mathcal{F}_{k-1}]  & \leq &   \frac{ \mathbf{E}[ \sigma_k^2\mathbf{1}_{\{|X_k| \geq v_k\}}|\mathcal{F}_{k-1}]+ \mathbf{E}[ X_k^2\mathbf{1}_{\{|X_k| \geq v_k\}}|\mathcal{F}_{k-1}]}{v_k}
+2\frac{ \mathbf{E}[(\sigma_k^2 |X_k|+|X_k|^3)\mathbf{1}_{\{|X_k| <v_k\}}  |\mathcal{F}_{k-1}]}{v_k^2}  \nonumber \\
    &=& \frac{  \mathbf{E}[(\sigma_k^2 +X_k^2)\mathbf{1}_{\{|X_k| \geq v_k\}}|\mathcal{F}_{k-1}]}{v_k}  +2\frac{ \mathbf{E}[(\sigma_k^2 |X_k|+|X_k|^3)\mathbf{1}_{\{|X_k| <v_k\}}  |\mathcal{F}_{k-1}]}{v_k^2} .\label{fsdfg}
	\end{eqnarray}
Thus, we have
	\begin{eqnarray*}
 \textbf{W}\big(\mathcal{L}(S_n), \mathcal{L}( N(0, s_n^2) )\big)  &\leq& |\mathbf{E}[h(S_n)-h(Z)]| \\
    & \leq&   \Big|\mathbf{E}[h(S_n+aZ')-h(Z+aZ')] \Big|+2a\\
	&  =&  \Big|\mathbf{E} \sum_{k=1}^n \mathbf{E}[R_k|\mathcal{F}_{k-1}] \Big|+2a\\
	& \leq&  \sum_{k=1}^n \bigg(\mathbf{E}\bigg[\frac{(\sigma_k^2  + X_k^2)\mathbf{1}_{\{|X_k| \geq v_k\}} }{v_k}\bigg]   +2 \mathbf{E}\bigg[\frac{(\sigma_k^2 |X_k|+|X_k|^3)\mathbf{1}_{\{|X_k| <v_k\}} }{v_k^2}\bigg]\bigg) + 2a.
	\end{eqnarray*}
Normalized $S_n$, we get
	\begin{eqnarray*}
 \textbf{W}\big(  S_n/s_n  \big)   \leq
	\frac{1}{s_n} \sum_{k=1}^n  \bigg(\mathbf{E}\bigg[\frac{(\sigma_k^2  + X_k^2)\mathbf{1}_{\{|X_k| \geq v_k\}} }{v_k}\bigg]   +2 \mathbf{E}\bigg[\frac{(\sigma_k^2 |X_k|+|X_k|^3)\mathbf{1}_{\{|X_k| <v_k\}} }{v_k^2}\bigg]\bigg)      +\frac{2a}{s_n},
	\end{eqnarray*}
which gives the desired inequality.

\subsection{Proof of Corollary \ref{fdgs}}
By  H\"{o}lder's inequality and Jensen's inequality, we get
\begin{eqnarray}
\mathbf{E}\bigg[\frac{ \sigma_k^2   \mathbf{1}_{\{|X_k| \geq v_k\}} }{v_k}\bigg] \!\!&\leq&\!\! \mathbf{E} \frac{ \sigma_k^2  |X_k|^{\delta}  }{v_k v_k\,\!\!^{\delta}} =  \mathbf{E} \frac{ \sigma_k^2  |X_k|^{\delta}  }{(\rho_k^2+a^2)^{(1+\delta)/2} } \nonumber\\
&\leq&\!\! \Big[\mathbf{E} \Big(\frac{ \sigma_k^2  }{\ (\rho_k^2+a^2)^{(1+\delta)/(2+\delta)}  }\Big)^{(2+\delta)/2} \Big]^{2/(2+\delta)} \Big[\mathbf{E}\Big( \frac{  |X_k|^{\delta}  }{  (\rho_k^2+a^2)^{ }\,\!\!^{\delta(1+\delta)/(4+2\delta)}} \Big)^{(2+\delta)/\delta}\Big]^{\delta/(2+\delta)} \nonumber \\
&\leq&\!\!\bigg[\mathbf{E} \frac{ |X_k|^{2+ \delta} }{ (\rho_k^2+a^2 )^{(1+ \delta)/2} } \bigg]^{2/(2+\delta)} \bigg[\mathbf{E}\frac{ |X_k|^{2+ \delta} }{ (\rho_k^2+a^2 )^{(1+ \delta)/2} }\bigg]^{\delta/(2+\delta)} \nonumber\\
&\leq&\!\! \mathbf{E} \frac{ |X_k|^{2+ \delta} }{ (\rho_k^2+a^2 )^{(1+ \delta)/2} }  . \label{ines21}
\end{eqnarray}
Similarly, we can prove that
\begin{eqnarray}
  \mathbf{E}\bigg[\frac{ \sigma_k^2 |X_k| \mathbf{1}_{\{|X_k| <v_k\}} }{v_k^2}\bigg] &\leq& \mathbf{E} \frac{ \sigma_k^2  |X_k|^{\delta}   v_k\,\!\!^{1-\delta}}{ \rho_k^2+a^2 } = \mathbf{E} \frac{ \sigma_k^2  |X_k|^{\delta}  }{ (\rho_k^2+a^2)^{(1+\delta)/2} } \nonumber \\
&\leq& \mathbf{E} \frac{ |X_k|^{2+ \delta} }{ (\rho_k^2+a^2 )^{(1+ \delta)/2} }  .
\end{eqnarray}
The following inequalities hold obviously
\begin{eqnarray}
 \mathbf{E}\bigg[\frac{  X_k^2 \mathbf{1}_{\{|X_k| \geq v_k\}} }{v_k}\bigg] \leq  \mathbf{E} \frac{ |X_k|^{2+\delta} }{v_k  v_k\,\!\!^{\delta}} =  \mathbf{E} \frac{ |X_k|^{2+ \delta} }{ (\rho_k^2+a^2 )^{(1+ \delta)/2} },\ \ \
	\end{eqnarray}
\begin{eqnarray}
	 \mathbf{E}\bigg[\frac{ |X_k|^3 \mathbf{1}_{\{|X_k| < v_k\}} }{v_k^2}\bigg] \leq  \mathbf{E} \frac{ |X_k|^{2+\delta} v_k\,\!^{1-\delta}   }{v_k^2} =  \mathbf{E} \frac{ |X_k|^{2+ \delta} }{ (\rho_k^2+a^2 )^{(1+ \rho)/2} }  . \label{ines25}
	\end{eqnarray}
Applying the inequalities (\ref{ines21}) - (\ref{ines25}) to Theorem \ref{theorem2.1}, we obtain the first desired inequality:
\begin{eqnarray*}
 \textbf{W} \big(  S_n/s_n  \big) & \leq& \frac{6}{s_n} \sum_{k=1}^n    \mathbf{E} \frac{ |X_k|^{2+\delta} }{  (\rho_k^2+a^2 )^{(1+ \delta)/2} }        +\frac{ 2a}{s_n}.
\end{eqnarray*}
This completes the proof of Corollary \ref{fdgs}.

\subsection{Proof of Corollary \ref{fdfff}}
By Jensen's inequality,   $\mathbf{E}[ |X_k|^{2+\delta} | \mathcal{F}_{k-1} ] \leq \gamma \, \mathbf{E}[ X_k^2 | \mathcal{F}_{k-1} ]$ and $\mathbf{E}[ X_k^2 | \mathcal{F}_{k-1} ] = \sigma_k^2$ together implies that for all $1 \leq k \leq n,$
$$ (\sigma_k^2)^{(2+\delta) /2 }  \leq  \mathbf{E}[ |X_k|^{2+\delta} | \mathcal{F}_{k-1} ] \leq  \gamma \, \sigma_k^2 ,$$
which  implies that $\sigma_k^2 \leq \gamma^{2/\delta}.$
Thus  $\sigma_k^2/s_n^2 \rightarrow 0, n \rightarrow \infty,$ uniformly for all $1 \leq k \leq n.$
Notice that $ \rho_{k }^2  = \sum_{i=k }^n \sigma_i^2.$ Taking  $a=1,$
by the conditions in Corollary \ref{fdfff},  it holds for $\delta \in (0, 1)$,
\begin{eqnarray*}
   \sum_{k=1}^n    \mathbf{E} \frac{ |X_k|^{2+\delta} }{ (\rho_k^2+1)^{(1+ \delta)/2} }     & =&  \mathbf{E} \sum_{k=1}^n   \frac{  \mathbf{E}[ |X_k|^{2+\delta} | \mathcal{F}_{k-1}   ] }{ (\rho_k^2  +1 )^{(1+ \delta)/2}}  \\
    & \leq &  s_n^{1-\delta} \gamma   \  \mathbf{E}\sum_{k=1}^n    \frac{  \sigma_k^2/s_n^2 }{ (\rho_k^2/s_n^2 +1/s_n^2 )^{(1+ \delta)/2}}      \\
     & \leq &  s_n^{1-\delta} \gamma\,  \Big(  \int_{0}^1\frac{1}{x^{(1+ \delta)/2}}dx + \frac{\gamma^2}{s_n^2} \Big)    \\
 &\leq& C_1\,s_n^{1-\delta} .
\end{eqnarray*}
Applying the last inequality to Corollary \ref{fdgs} with $a=1$,  it holds for $\delta \in (0, 1)$,
     \begin{eqnarray}
	 \textbf{W} \big( S_n/s_n  \big) \leq    C\, \frac{1}{ s_n^{\delta }}.
	\end{eqnarray}
When $\delta =1$,  we have
\begin{eqnarray*}
   \sum_{k=1}^n    \mathbf{E} \frac{ |X_k|^{3} }{  \rho_k^2+1  }     & =&  \mathbf{E} \sum_{k=1}^n    \frac{  \mathbf{E}[ |X_k|^3 | \mathcal{F}_{k-1}   ] }{  \rho_k^2  +1  } \,\leq \,   \gamma   \, \mathbf{E}\sum_{k=1}^n    \frac{  \sigma_k^2/s_n^2 }{  \rho_k^2/s_n^2 +1/s_n^2 }
  \, \leq \,  C_1\,  \ln s_n .
\end{eqnarray*}
Applying the last inequality to Corollary \ref{fdgs} and taking $a=1$, we get for $\delta =1$,
     \begin{eqnarray}
	 \textbf{W} \big( S_n/s_n  \big) \leq    C\, \frac{\ln s_n}{ s_n }.
	\end{eqnarray}
This completes the proof of Corollary \ref{fdfff}.

\subsection{Proof of Corollary \ref{fdffs}}
We only give a proof for the case $\delta \in (0, 1)$. For the case $\delta=1$, the proof is similar.
Taking $a=\alpha,$ by the conditions in Corollary \ref{fdffs},  we get
\begin{eqnarray*}
  6  \sum_{k=1}^n    \mathbf{E} \frac{ |X_k|^{2+\delta} }{ (\rho_k^2+\alpha )^{(1+ \delta)/2} }        + 2 \alpha & \leq&   \sum_{k=1}^n   \frac{ 6}{ ((n-k) \alpha +\alpha )^{(1+ \delta)/2} } \mathbf{E} |X_k|^{2+\delta}       +  2 \alpha    \\
 &\leq& C_1 \bigg( \sum_{k=1}^n   \frac{ 1}{ ((n-k)   + 1)^{(1+ \delta)/2} }       +  1  \bigg) \\
  &\leq& C_2 n^{(1-\delta)/2} .
\end{eqnarray*}
Again by the conditions in Corollary \ref{fdffs}, we have $s_n^2  \geq n\alpha$. Therefore, it holds
\begin{eqnarray*}
 \frac{6}{s_n} \sum_{k=1}^n    \mathbf{E} \frac{ |X_k|^{2+\delta} }{ (\rho_k^2+a^2 )^{(1+ \delta)/2} }        +\frac{ 2a}{s_n}& \leq& C_2 \frac{1}{\sqrt{n\alpha}} n^{(1-\delta)/2} \leq C_3 n^{ -\delta /2}.
\end{eqnarray*}
Applying the last inequality to Corollary \ref{fdgs}, we get the desired inequality.

\subsection{Proof of Theorem \ref{theorem3.2}}

The proof of Theorem \ref{theorem3.2} is similar to that of Theorem \ref{theorem2.1}, with $\sigma_k$ replaced by $\overline{\sigma}_k$.
Similar to the proof of (\ref{fsdfg}), we obtain
	\begin{eqnarray*}
	 \Big|\mathbf{E}[R_k|\mathcal{F}_{k-1}] \Big|
	&\leq& \Big|\mathbf{E}[(X_k^2-\overline{\sigma}_k^2) f'_{S_{k-1}, \tau_k}(S_{k-1} +T'_{k+1})|\mathcal{F}_{k-1}] \Big| \\
&&  + \frac{\overline{\sigma}_k^2 \mathbf{E}[ \mathbf{1}_{\{|X_k| \geq \tau_k\}}|\mathcal{F}_{k-1}]}{\sqrt{\overline{\rho}_k^2+a^2}} + \frac{ \mathbf{E}[ X_k^2\mathbf{1}_{\{|X_k| \geq \tau_k\}}|\mathcal{F}_{k-1}]}{\sqrt{\overline{\rho}_k^2+a^2}}
+2\frac{ \mathbf{E}[|X_k|^3\mathbf{1}_{\{|X_k| < \tau_k\}}  |\mathcal{F}_{k-1}]}{\overline{\rho}_k^2+a^2} \\
    &\leq& \frac{1}{\sqrt{\overline{\rho}_k^2+a^2}}  \Big|\mathbf{E}[ X_k^2   |\mathcal{F}_{k-1}]-\overline{\sigma}_k^2 \Big|   +  \frac{  \mathbf{E}[(\overline{\sigma}_k^2 +X_k^2)\mathbf{1}_{\{|X_k| \geq \tau_k\}}|\mathcal{F}_{k-1}]}{\sqrt{\overline{\rho}_k^2+a^2}} \\
     && +2\frac{ \mathbf{E}[(\overline{\sigma}_k^2 |X_k|+|X_k|^3)\mathbf{1}_{\{|X_k| < \tau_k\}}  |\mathcal{F}_{k-1}]}{\overline{\rho}_k^2+a^2}.
	\end{eqnarray*}
Then, it holds for any $a\geq0$,
	\begin{eqnarray*}
|\mathbf{E}[h(S_n)-h(Z)]|
    &\leq& |\mathbf{E}[h(S_n+aZ')-h(Z+aZ')]|+2a\\
	&  =& \ \Big|\mathbf{E} \sum_{k=1}^n \mathbf{E}[R_k|\mathcal{F}_{k-1}] \Big|+2a\\
	&  \leq& \ \sum_{k=1}^n \Bigg(\frac{ 1}{\sqrt{\overline{\rho}_k^2+a^2}}\bigg(\mathbf{E} \big|\mathbf{E}[ X_k^2   |\mathcal{F}_{k-1}]-\overline{\sigma}_k^2 \big| +\overline{\sigma}_k^2 \mathbf{E}[ \mathbf{1}_{\{|X_k| \geq \tau_k\}} ] +  \mathbf{E}[ X_k^2\mathbf{1}_{\{|X_k| \geq \tau_k\}}]\bigg) \\
& &  +2 \frac{\overline{\sigma}_k^2\mathbf{E}[  |X_k| \mathbf{1}_{\{|X_k| < \tau_k\}}  ]+\mathbf{E}[|X_k|^3\mathbf{1}_{\{|X_k| < \tau_k\}}  ]}{\overline{\rho}_k^2+a^2}\Bigg) + 2a,
	\end{eqnarray*}
where  $h$ can be any Lipschitz function such that $||h'||_{\infty} \leq 1$.
Thus, we have
	\begin{eqnarray*}
&& \textbf{W}\Big( \mathcal{L} (S_n),   N(0, s_n^2)  \Big) \\
 &&\leq  |\mathbf{E}[h(S_n)-h(Z)]| \leq  \ \Big|\mathbf{E}[h(S_n+aZ')-h(Z+aZ')] \Big|+2a\\
	& & \leq \ \sum_{k=1}^n \Bigg(\frac{ 1}{\sqrt{\overline{\rho}_k^2+a^2}}\bigg(\mathbf{E} \big|\mathbf{E}[ X_k^2   |\mathcal{F}_{k-1}]-\overline{\sigma}_k^2 \big| +\overline{\sigma}_k^2 \mathbf{E}[ \mathbf{1}_{\{|X_k| \geq \tau_k\}} ] +  \mathbf{E}[ X_k^2\mathbf{1}_{\{|X_k| \geq \tau_k\}}]\bigg) \\
& & \ \ \ \ \ + \ 2 \frac{\overline{\sigma}_k^2\mathbf{E}[  |X_k| \mathbf{1}_{\{|X_k| < \tau_k\}}  ]+\mathbf{E}[|X_k|^3\mathbf{1}_{\{|X_k| < \tau_k\}}  ]}{\overline{\rho}_k^2+a^2}\Bigg) + 2a.
	\end{eqnarray*}
Normalized $S_n$, we obtain
	\begin{eqnarray*}
\textbf{W} \big(  S_n/s_n   \big)  &\leq&   \frac{1}{s_n} \sum_{k=1}^n \Bigg(\frac{ 1}{\sqrt{\overline{\rho}_k^2+a^2}}\Big(\mathbf{E} \big|\mathbf{E}[ X_k^2   |\mathcal{F}_{k-1}]-\overline{\sigma}_k^2 \big| +\overline{\sigma}_k^2 \mathbf{E}[ \mathbf{1}_{\{|X_k| \geq \tau_k\}} ] +  \mathbf{E}[ X_k^2\mathbf{1}_{\{|X_k| \geq \tau_k\}}]\Big) \\
&&\  +\ 2 \frac{\overline{\sigma}_k^2\mathbf{E}[  |X_k| \mathbf{1}_{\{|X_k| < \tau_k\}}  ]+\mathbf{E}[|X_k|^3\mathbf{1}_{\{|X_k| < \tau_k\}}  ]}{\overline{\rho}_k^2+a^2}\Bigg) +\frac{2a}{s_n}\\
& \leq&   \frac{1}{s_n} \sum_{k=1}^n \Bigg(\frac{ 1}{  \tau_{ k}}\Big(\mathbf{E} \big|\mathbf{E}[ X_k^2   |\mathcal{F}_{k-1}]-\overline{\sigma}_k^2 \big| +  \mathbf{E}[ X_k^2\mathbf{1}_{\{|X_k| \geq \tau_k\}}]\Big) \\
&&  \   + \ \frac{ 1}{ \tau_{k}^2 }\Big(3\,\overline{\sigma}_k^2  \mathbf{E}  |X_k| + 2\,\mathbf{E}[|X_k|^3\mathbf{1}_{\{|X_k| < \tau_k\}}  ] \Big)\Bigg) +\frac{2a}{s_n},
	\end{eqnarray*}
which gives the first desired inequality.

\subsection{Proof of Corollary \ref{fsvg}}

By the fact $\mathbf{E} \big|\mathbf{E}[ X_k^2   |\mathcal{F}_{k-1}]-\overline{\sigma}_k^2 \big|=O( k^{-\alpha}) \   (k\rightarrow \infty)$ and $s_n^2-s_k^2  \asymp n- k$ for any $n >k$, we can deduce that
  \begin{eqnarray*}
  \sum_{k=1}^n  \frac{ \mathbf{E} \big|\mathbf{E}[ X_k^2   |\mathcal{F}_{k-1}]-\overline{\sigma}_k^2 \big|  }{ \tau_k  }   &\leq& C\, \sum_{k=1}^n  \frac{ 1 }{ \sqrt{ s_n^2 -s_{k-1}^2 + a^2 } \ } k^{-\alpha}  \\
  &\leq& C\, \sum_{k=1}^n  \frac{ 1 }{ \sqrt{ n -k + 1 } \ } k^{-\alpha} \\
  &=& C\, \sum_{k=1}^n  \frac{ 1 }{ \sqrt{ 1 -k/n + 1/n } \ } \frac{1}{ (k/n)^{ \delta/2}} \, \frac{1}{n} \, n^{(1-2\alpha )/2} \\
  &\leq& C\, \int_{0}^1 \frac{}{}\frac{}{}   \frac{ 1 }{ \sqrt{ 1 -x   } \ } \frac{1}{ x^{ \delta/2}} \,dx \, n^{(1-2\alpha )/2} \\
  &\leq& C\,  n^{(1-2\alpha)/2}.
\end{eqnarray*}
Applying the last inequality to  (\ref{gsds}), we get the first desired inequality.
Taking $a=1,$ by the condition in Corollary \ref{fsvg}, it holds for $\delta \in (0, 1)$,
\begin{eqnarray*}
   \frac{1}{ \sqrt{n}} \sum_{k=1}^n \frac{  \mathbf{E}[ X_k^2\mathbf{1}_{\{|X_k| \geq \tau_k\}}] }{ \tau_k  }  \, & \leq &  \frac{1}{ \sqrt{n}}\sum_{k=1}^n    \frac{\mathbf{E} |X_k|^{2+\delta} }{ (\overline{\rho}_k^2+1)^{(1+ \delta)/2} }      \\
    & \leq & \frac{  C}{ \sqrt{n}}\sum_{k=1}^n      \frac{  1 }{ (n-k+1 )^{(1+ \delta)/2}}      \\
     & \leq &  \frac{  C}{  n^{\delta/2}}.
\end{eqnarray*}
Similarly, it holds  for  $\delta \in (0, 1)$,
\begin{eqnarray*}
 \frac{1}{ \sqrt{n}} \sum_{k=1}^n   \, \frac{ \mathbf{E}[|X_k|^3\mathbf{1}_{\{|X_k| < \tau_k\}}  ]  }{  \tau_k\,\!^2 }     & \leq &  \frac{1}{ \sqrt{n}}\sum_{k=1}^n    \frac{\mathbf{E} |X_k|^{2+\delta} }{ (\overline{\rho}_k^2+1)^{(1+ \delta)/2} }
     \ \leq \  \frac{  C}{  n^{\delta/2}}.
\end{eqnarray*}
Applying the foregoing results to Corollary \ref{fsvg},  we get for $\delta \in (0, 1)$,
     \begin{eqnarray}\nonumber
	 \textbf{W} \big( S_n/s_n \big) \leq  \frac{  C}{  n^{\delta/2}}.
	\end{eqnarray}
For $\delta =1$, Corollary \ref{fsvg} can be proved similarly.

\subsection{Proof of Theorem \ref{th00s1}}

Denote   $c_\delta$ a positive constant depending  on   $\rho,  \upsilon^2 , \sigma^2, \tau^2,$  $\mu$, $ \nu, \alpha,$  $\mathbf{E} V^{-\alpha},$ $\mathbf{E}| Z_1  -m_0   |^{2+\delta }$ and $\mathbf{E}   |  m_0 -m   |^{2+\delta }$. In particular,   it does not depend on $n$ and $x$. 
Moreover, the exact values of $c_{\delta}$ may vary from line to line.
\begin{lemma}  \label{thsn20}
Assume the conditions of Theorem \ref{th00s1}.
   It holds for all $0<x \leq \mu \sqrt{n} /\nu,$
$$\mathbf{P}\bigg( \Big| \frac{\ln Z_{n}  - \mu  n }{\sqrt{ n}\, \nu }\Big| \geq x \bigg) \leq  c_{\delta}^{-1} \exp\bigg\{    - c_{\delta}\, x^2 \bigg\}.$$
\end{lemma}
\emph{Proof}. It is easy to see that for all  $x> 0,$
\begin{eqnarray}\label{dfsafh}
\mathbf{P}\bigg( \Big| \frac{\ln Z_{n}  - \mu  n }{\sqrt{ n}\, \nu }\Big| \geq x \bigg)  =  I_1 +I_2,
\end{eqnarray}
where
$$I_1=  \mathbf{P}\bigg(   \frac{\ln Z_{n}  - \mu  n }{\sqrt{ n}\, \nu }  \leq - x \bigg) \ \ \ \   \textrm{and}  \ \ \ \   I_2 = \mathbf{P}\bigg(   \frac{\ln Z_{n}  - \mu  n }{\sqrt{ n}\, \nu }  \geq x \bigg).$$
Clearly, we have  the following decomposition:
\begin{equation}\label{decopos}
 \ln Z_n = \sum_{i=1}^n X_i + \ln V_n,
\end{equation}
where $X_i=\ln m_{i-1} (i\geq 1)$ are i.i.d.\ random variables depending only on the environment $\xi.$
Denote $\eta_{i}=(X_i-\mu)/\sqrt{n}\, \nu.$
Clearly, it holds for all  $  x >0,$
 \begin{eqnarray}
   I_1
  &=&    \mathbf{P}\bigg(  \sum_{i=1}^n \eta_{i} + \frac{\ln V_{n}}{  \sqrt{n}\, \nu} \leq - x \bigg )   \nonumber \\
   &\leq&   \mathbf{P}\bigg( \sum_{i=1}^n \eta_{i}  \leq - \frac x 2 \bigg ) + \mathbf{P}\bigg(  \frac{\ln V_{n}}{  \sqrt{n}\, \nu} \leq -  \frac x 2  \bigg ) . \label{dasa01}
\end{eqnarray}
The condition $\mathbf{E}   |  m_0 -m   |^{2+\rho } < \infty$ implies that $  \mathbf{E} e^{c_\rho\, |X_1| }  < \infty$ for
some $c_\rho > 0$.
By Corollary 1.2 of Liu and Watbled \cite{LW09}, it holds
\begin{displaymath}
  \mathbf{P}\bigg( \frac{1}{n}\sum_{i=1}^n X_i  -\mu \leq - x \bigg )  \leq \left\{ \begin{array}{ll}
\exp\Big\{ - c_\rho \,  n \, x^2    \Big\},\ \ \ \  & \textrm{ $0< x\leq1$,}\\
  &  \\
\exp\Big\{ - c_\rho \, n \, x    \Big\},\ \ \ \  & \textrm{ $x> 1$. }
\end{array} \right.
\end{displaymath}
The last inequality implies that  for all $ 0  <x \leq \mu \sqrt{n} /\nu, $
 \begin{eqnarray}
 \mathbf{P}\bigg( \sum_{i=1}^n \eta_{i}  \leq - \frac x 2 \bigg )
 \, \leq \, \exp\bigg\{ - c_\rho \,  x^2    \bigg\}.
\label{dasa02}
\end{eqnarray}
By Markov's inequality,    it is easy to see that for  all $x>0,$
 \begin{eqnarray}
\mathbf{P}\bigg(  \frac{\ln V_{n}}{  \sqrt{n}\, \nu} \leq -  \frac x 2  \bigg )   = \mathbf{P}\bigg(  V_{n}^{-\alpha  } \geq \exp\Big\{\frac x 2 \sqrt{n} \nu \alpha \Big\} \bigg)
 \leq   \exp\bigg\{ - \frac x 2 \sqrt{n} \nu \alpha  \bigg \}  \mathbf{E}  V_{ n}  ^{-\alpha } ,\nonumber
\end{eqnarray}
where $\alpha$ is given by (\ref{defv55}).
As $V_n \rightarrow V$ in $\mathbf{L}^{1} ,$   we have $V_n=\mathbf{E}[V| \mathcal{F}_n]$ a.s.\ Using Jensen's inequality,
we can deduce that
 \begin{eqnarray*}
 V_{ n} ^{-\alpha } =(\mathbf{E}[V  | \mathcal{F}_{ n}])^{-\alpha  } \leq \mathbf{E}[  V  ^{-\alpha } | \mathcal{F}_{ n}].
\end{eqnarray*}
Taking expectations on both sides of the last inequality, by (\ref{defv55}),
we get $\mathbf{E} V_{n} ^{-\alpha} \leq \mathbf{E}  V ^{-\alpha  }  < \infty.$
Hence, we have for all  $ 0  <x \leq \mu \sqrt{n} /\nu, $
 \begin{eqnarray}
\mathbf{P}\bigg(  \frac{\ln V_{n}}{  \sqrt{n}\, \nu} \leq -  \frac x 2  \bigg )
 \leq     c_\rho \, \exp\bigg\{ - \frac x 2 \sqrt{n} \nu \alpha  \bigg \}  \leq  c_\rho \,  \exp\bigg\{ - c_\rho '  x^2    \bigg\} .   \label{dasa03}
\end{eqnarray}
Combining (\ref{dasa01})-(\ref{dasa03}) together, we get   for  all $ 0  <x \leq \mu \sqrt{n} /\nu, $
 \begin{eqnarray}
I_1   \leq    c_\rho^{-1}  \exp\Big\{ - c_\rho    x^2    \Big\} .\nonumber
\end{eqnarray}
Similarly, we can prove that the last bound holds also for $I_2.$
The desired inequality follows by (\ref{dfsafh}) and
the upper bounds of $I_1$ and $I_2$.  This completes the proof of Lemma \ref{thsn20}.

Now, we are in position to prove Theorem \ref{th00s1}.
Denote
 \begin{eqnarray}\label{gdfbg}
  \tilde{\xi}_{k+1}=  \frac{Z_{ k+1}}{Z_{ k}}  -m  ,
 \end{eqnarray}
$\mathfrak{F}_{n_0} =\{ \emptyset, \Omega \}  $ and $\mathfrak{F}_{k+1}=\sigma \{ \xi_{i-1}, Z_{i}: n_0\leq i\leq k+1  \}$ for  all $k\geq n_0$.
Notice that $X_{k,i}$ is independent of $Z_k.$
Then it is easy to verify that
 \begin{eqnarray}\nonumber
\mathbf{E}[ \tilde{\xi}_{k+1}  |\mathfrak{F}_{k } ]  =   Z_{ k} ^{-1 }  \mathbf{E}[   Z_{ k+1}  -mZ_{ k}  |\mathfrak{F}_{k } ] \ = \  Z_{ k} ^{-1 } \sum_{i=1}^{Z_k} \mathbf{E}[   X_{ k, i}  -m   |\mathfrak{F}_{k } ]    =    Z_{ k} ^{-1 } \sum_{i=1}^{Z_k} \mathbf{E}[   X_{ k, i}  -m   ] \  = \ 0.
\end{eqnarray}
 Thus
 $(\tilde{\xi}_k, \mathfrak{F}_k)_{k=n_0+1,...,n_0+n}$ is a finite sequence of martingale  differences.
 Moreover, it is easy to see that
\begin{eqnarray}
   \mathbf{E}[ \tilde{\xi}_{k+1}^2  |\mathfrak{F}_{k } ] &=&
   Z_{k}^{-2} \mathbf{E}[  ( Z_{ k+1} -mZ_{ k} )^2  |\mathfrak{F}_{k } ]=   Z_{k}^{-2} \mathbf{E}\Big[  \Big( \sum_{i=1}^{Z_k}    (X_{ k, i}  -m)  \Big)^2  \Big|\mathfrak{F}_{k } \Big]\nonumber \\
&=&   Z_{k}^{-2} \mathbf{E}\Big[  \Big( \sum_{i=1}^{Z_k}    (X_{ k, i} -m_k+ m_k -m)  \Big)^2  \Big|\mathfrak{F}_{k } \Big]  \nonumber  \\
&=&  Z_{k}^{-2} \Bigg(\mathbf{E}\Big[  \Big( \sum_{i=1}^{Z_k}    (X_{ k, i} -m_k)    \Big)^2 \Big |\mathfrak{F}_{k } \Big ] +\, \mathbf{E}\Big[  \Big( \sum_{i=1}^{Z_k}    (m_k -m ) \Big)^2  \Big|\mathfrak{F}_{k } \Big] \nonumber\\
  &&  +\, 2\  \mathbf{E}\Big[  \Big( \sum_{i=1}^{Z_k}    (X_{ k, i} -m_k)    \Big)\Big( \sum_{i=1}^{Z_k}    (m_k -m ) \Big)  \Big|\mathfrak{F}_{k } \Big]\Bigg).\nonumber
\end{eqnarray}
Notice that conditionally
on  $\xi_n$, the random variables $(X_{n,i})_{i\geq1}$ are i.i.d. Thus,  we can deduce that
\begin{eqnarray*}
\mathbf{E}\Big[  \Big( \sum_{i=1}^{Z_k}    (X_{ k, i} -m_k)    \Big)^2  \Big |\mathfrak{F}_{k }  \Big]&=& \mathbf{E}\Big[ \mathbf{E}\Big[  \Big( \sum_{i=1}^{Z_k}    (X_{ k, i} -m_k)    \Big)^2  \Big|\xi_k , \mathfrak{F}_{k } \Big]   \Big| \mathfrak{F}_{k } \Big] \\
&=& \mathbf{E}\Big[    \sum_{i=1}^{Z_k}   \mathbf{E}[ (X_{ k, i} -m_k)^2      | \xi_k , \mathfrak{F}_{k } ]   \Big| \mathfrak{F}_{k } \Big] \ = \ \sum_{i=1}^{Z_k}    \mathbf{E}[ (X_{ k, i} -m_k)^2      |   \mathfrak{F}_{k } ]\\
&=&  \sum_{i=1}^{Z_k}   \mathbf{E}  (X_{ k, i} -m_k)^2 \ = \ Z_k  \tau^2.
\end{eqnarray*}
Similarly, we can prove that
\begin{eqnarray*}
\mathbf{E}\Big[  \Big( \sum_{i=1}^{Z_k}    (X_{ k, i} -m_k)    \Big)\Big( \sum_{i=1}^{Z_k}    (m_k -m ) \Big)  \Big| \mathfrak{F}_{k } \Big]
 =  0
\end{eqnarray*}
and
\begin{eqnarray*}
\mathbf{E}\Big[  \Big( \sum_{i=1}^{Z_k}    (m_k -m ) \Big)^2  \Big|\mathfrak{F}_{k } \Big] =  Z_k^2 \mathbf{E}     (m_k -m ) ^2= Z_k^2 \mathbf{E}     (m_0 -m ) ^2=  Z_k^2\sigma^2 .
\end{eqnarray*}
Therefore, it holds
\begin{eqnarray}
  \mathbf{E}[ \tilde{\xi}_{k+1}^2  |\mathfrak{F}_{k } ] = Z_k^{-1}  \tau^2 +  \sigma^2 \geq   \sigma^2 .
\end{eqnarray}
By the last line and the fact $Z_k \geq 1$ a.s.\ for all $k\geq1$, it is easy to see that
\begin{eqnarray}
 \mathbf{E}  \big|\mathbf{E}[  \tilde{\xi}_{k+1} ^{2 } |\mathfrak{F}_{k } ] -  \mathbf{E}  \tilde{\xi}_{k+1}^2 \big |    &\leq&\, 2\,\tau^2 \mathbf{E} Z_k^{-1} \ \leq \ 2\,\tau^2 \Big( \mathbf{P}( Z_k \leq k   )+ \frac1k \mathbf{P}( Z_k > k   ) \Big) \nonumber \\
 &\leq&\, 2\,\tau^2 \bigg( \mathbf{P}( Z_k \leq k   )+ \frac1k  \bigg).  \label{fsdsl3}
\end{eqnarray}
By Lemma \ref{thsn20}, we have
\begin{eqnarray}\label{fsdsl4}
 \mathbf{P}( Z_k \leq k   )  \leq  \mathbf{P}\bigg( \frac{\ln Z_{k}  - \mu  k }{\sqrt{ k}\, \nu }    \leq \frac{\ln k  - \mu  k }{\sqrt{ k}\, \nu }    \bigg)  \leq  c_{\rho}^{-1} \exp\Big\{    - c_{\rho}\, k \Big\}.
\end{eqnarray}
Applying (\ref{fsdsl4}) to (\ref{fsdsl3}), we get
\begin{eqnarray}
 \mathbf{E}  |\mathbf{E}[  \tilde{\xi}_{k+1} ^{2 } |\mathfrak{F}_{k } ] -  \mathbf{E}  \tilde{\xi}_{k+1}^2 |    \ll   \frac1k.
\end{eqnarray}
 Notice that
\begin{eqnarray}
  \mathbf{E}[ |\tilde{\xi}_{k+1}|^{2+\delta} |\mathfrak{F}_{k } ] &=&
  Z_{k}^{-2-\delta } \mathbf{E}[  | Z_{ k+1} -mZ_{ k} |^{2+\delta }  |\mathfrak{F}_{k } ]  \nonumber \\
& = &  Z_{k}^{-2-\delta }  \mathbf{E}\bigg[  \Big| \sum_{i=1}^{Z_k}    (X_{ k, i}  -m)  \Big|^{2+\delta }  \bigg|\mathfrak{F}_{k } \bigg]. \label{ines2.8}
\end{eqnarray}
By Minkowski's inequality, we have
\begin{eqnarray}
&& \mathbf{E}\bigg[  \Big| \sum_{i=1}^{Z_k}    (X_{ k, i}  -m)  \Big|^{2+\delta }  \bigg|\mathfrak{F}_{k } \bigg]
\leq  \Bigg(  \sum_{i=1}^{Z_k}  \Big( \mathbf{E}[  |  X_{ k, i}  -m   |^{2+\delta }  |\mathfrak{F}_{k } ]   \Big)^{\frac{1}{2+\delta} }  \Bigg)^{2+\delta} \nonumber \\
&&=\Bigg(  \sum_{i=1}^{Z_k}  \Big( \mathbf{E}[  |  X_{ k, i}  -m_k + m_k -m  |^{2+\delta }  |\mathfrak{F}_{k } ]   \Big)^{\frac{1}{2+\delta} }  \Bigg)^{2+\delta} \nonumber\\
&&\leq Z_k^{2+\delta}  2^{1+\delta}    \bigg( \mathbf{E}[  |  X_{ k, i}  -m_k   |^{2+\delta }  |\mathfrak{F}_{k } ] + \mathbf{E}[  |   m_k -m  |^{2+\delta }  |\mathfrak{F}_{k } ]  \bigg)   \nonumber \\
&& =   2^{1+\delta}  Z_k^{2+\delta}  \Big(      \mathbf{E}| Z_1  -m _0  |^{2+\delta } +  \mathbf{E}   |  m_0 -m   |^{2+\delta } \Big) . \label{dfsds4}
\end{eqnarray}
By (\ref{dfsds4}) and (\ref{ines2.8}), we get
 $$ \sup_{k\geq n_0} \Big\|  \mathbf{E}[ |\tilde{\xi}_{k+1}|^{2+\delta} |\mathfrak{F}_{k } ] \Big\|_\infty\leq 2^{1+\delta}   \Big(      \mathbf{E}| Z_1  -m _0  |^{2+\delta } +  \mathbf{E}   |  m_0 -m   |^{2+\delta } \Big) < \infty.$$
Applying Corollaries  \ref{fsdvs} and \ref{fsvg} to $(\tilde{\xi}_k, \mathfrak{F}_k)_{k=n_0+1,...,n_0+n}$, we obtain the desired results.
This completes the proof of Theorem \ref{th00s1}.


%
\section*{Acknowledgements}
%
The authors would like to thank Quansheng Liu for his helpful discussion on the harmonic moments for branching processes in a random environment.
Fan was partially supported by the National Natural Science Foundation
of China (Grant Nos.\,11971063 and 12371155). Su  was partially supported by the National Natural Science Foundation
of China (Grant Nos.\,12271457 and 11871425).

\end{document}